\documentclass[12pt]{article}
\usepackage{float,amsfonts,amssymb,makeidx}
\usepackage{amsmath}
\newtheorem{thm}{Theorem}[section]

\newtheorem{cor}[thm]{Corollary}
\newtheorem{prop}[thm]{Proposition}

\newtheorem{defin}[thm]{Definition}

\newbox\sample
\newif\ifproofmode
\newif\ifsymindex
\global\symindexfalse
\newwrite\inx
\def\indsyma#1#2{\ifproofmode\marginpar{$\scriptstyle#1$}\fi%
\ifx#2\empty\write\inx{$\noexpand#1$,\space\thepage}%
\write\inx{\string\newline}\else%
\write\inx{$\noexpand#1$,\space#2,\space\thepage}%
\write\inx{\string\newline}\fi\ignorespaces}%
\def\indsym#1#2{\ifsymindex%
\ifproofmode\marginpar{$\scriptstyle#1$}\fi%
\ifx#2\empty\write\inx{\string\item \space$\noexpand#1$,\space\thepage}%
\else%
\write\inx{\string\item \space$\noexpand#1$,\space#2,\space\thepage}%
\fi\ignorespaces\fi}%
\newskip\dangerskipb
\newskip\dangerskip
\def\hang{\hangindent\dangerskip}

\def\s#1{{\cal #1}}

\def\remark{\bigskip\noindent{\bf Remark:}\enspace}

\font\manual=manfnt at 12pt
\def\danbend{{\manual\char127}}
\def\datanger{\medbreak\begingroup\clubpenalty=10000
 \def\par{\endgraf\endgroup\medbreak} \noindent\hang\hangafter=-2
 \hbox to0pt{\hskip-3.5pc\danbend\hfill}}
\outer\def\danger{\datanger}%
\def\ddatanger{\medbreak\begingroup\clubpenalty=10000
 \def\par{\endgraf\endgroup\medbreak} \noindent\hang\hangafter=-2
 \hbox to0pt{\hskip-3.5pc\danbend\kern1pt%
\danbend\hfill}}

\def\dobdownarrow{\mathop{\vbox{\kern2pt \hbox{$\Big\downarrow$}\kern-16.5pt
                          \nointerlineskip\hbox{$\Big\downarrow$}}}}

\def\lrightarrow{\hbox to 25pt{\rightarrowfill}}

\def\supexp{exp(m,n,p)=m^{m^{m^{\cdot^{\cdot^{\cdot^{m^{p}}}}}}}
\vbox{\hbox{$\Big\}\scriptstyle n$}\kern0pt}}

\def\supexpo#1#2#3{#1^{#1^{\cdot^{\cdot^{\cdot^{#1^{#2}}}}}}
\vbox{\hbox{$\Big\}\scriptstyle #3$}\kern0pt}}

\def\sqr#1#2{{\vcenter{\hrule height .#2pt
         \hbox{\vrule width.#2pt height#1pt \kern#1pt
             \vrule width.#2pt}
         \hrule height.#2pt}}}

\def\co{\colon}

\newskip\bogcentering \bogcentering= 0pt plus 1000pt minus 1000pt 

\def\matth{\mathsurround=0pt}
\def\fakrightarrowfill{$\matth \mathord- \mkern-6mu
  \cleaders\hbox{$\mkern-2mu \mathord- \mkern-2mu$}\hfill
 \mkern-6mu \mathord\rightarrow$}

\def\fakoverrightarrow#1{\vbox{\ialign{##\crcr
  \fakrightarrowfill\crcr\noalign{\kern-1pt\nointerlineskip}
 $\hfil\displaystyle{#1}\hfil$\crcr}}}

\newif\ifdtatp

\def\displaty{%
\global \dtatptrue \openup \jot \matth \everycr{\noalign{\ifdtatp \global 
\dtatpfalse \vskip -\lineskiplimit \vskip \normallineskiplimit \else 
\penalty \interdisplaylinepenalty \fi }}}

\def\displaylignes#1{\displaty
   \halign{\hbox to\displaywidth{$\displaystyle##$}\crcr
   #1\crcr}}
\def\leqaligneno#1{\displaty \tabskip=\bogcentering
 \halign to\displaywidth{\hfil$\displaystyle{##}$\tabskip=0pt
 &$\displaystyle{{}##}$\hfil\tabskip=\bogcentering
 &\kern-\displaywidth\rlap{$##$}\tabskip=\displaywidthpt\crcr
 #1\crcr}}
\def\ligne{\hbox to\hsize}
\newdimen\nouvpagewidth
\newdimen\offwidth
\newdimen\lawidthoui
\nouvpagewidth=\lawidthoui
\def\kboxit#1{\vbox{\hrule\hbox{\vrule\kern3pt
              \vbox{\kern3pt#1\kern3pt}\kern3pt\vrule}\hrule}}

\def\kboxitb#1{\vbox{\hrule\hbox{\vrule\kern3pt
              \vbox{\kern3pt#1\kern3pt}\kern3pt\vrule}\hrule}}

\def\laboxaround#1{
\aboxaround{\hbox to\hsize{\hfill\box2\hfill}}{#1}
}

\def\boxar#1#2{
\aboxaround{\hbox to\hsize{\hfill#1\hfill}}{#2}
}

\def\aboxaround#1#2{
\setbox4=\vbox{\hsize #2\noindent\strut#1\strut}
\kboxitb{\box4}}

\def\kframeit#1{\vbox{\hrule\hbox{\vrule\kern5pt
              \vbox{\kern5pt#1\kern5pt}\kern5pt\vrule}\hrule}}

\newskip\savnormalbaselineskip
\newskip\savnormallineskip
\newdimen\savnormallineskiplimit

\def\reals{\mathbb{R}}
\def\complex{\mathbb{C}}
\def\integs{\mathbb{Z}}
\let\Immag=\Im
\def\co{\colon}
\def\mapdef#1#2#3{#1\co #2\rightarrow #3}

\def\s#1{{\cal #1}}
\def\norme#1{\left\|#1\right\|}
\def\remark{\bigskip\noindent{\bf Remark:}\enspace}

\def\id{\mathrm{id}}
\def\amsmata#1#2#3#4{%
\begin{pmatrix}
#1 & #2\\
#3 & #4
\end{pmatrix}
}

\begin{document}
\title{Logarithms and Square Roots of Real Matrices\\
Existence, Uniqueness, and  Applications\\
in Medical Imaging}
\author{Jean Gallier \\
%
 \\
Department of Computer and Information Science\\
University of Pennsylvania\\
Philadelphia, PA 19104, USA\\
{\tt jean@cis.upenn.edu
}
}
\maketitle
\vspace{0.3cm}
\noindent
{\bf Abstract.}
The need for computing logarithms or square roots of
real matrices arises in a number of applied problems.
A significant class of  problems comes from  medical imaging.
One of these problems is
to interpolate and to perform statistics on data represented
by certain kinds of matrices (such as symmetric positive definite
matrices in DTI). Another important and difficult problem
is the registration of medical images. 
For both of these problems, the ability to compute
logarithms of real matrices turns out to be crucial.
However, not all real matrices have a {\it real\/}
logarithm and thus, it is important to have sufficient
conditions for the existence (and possibly the uniqueness)
of a real logarithm for a real matrix. 
Such conditions (involving the eigenvalues
of a matrix) are known, both for the logarithm and the square
root.

\medskip
As far as I know, with the exception of Higham's recent book
\cite{HighamFM},
proofs of the results involving these conditions
are scattered in the literature
and it is not easy to locate them. Moreover, Higham's excellent book
assumes a certain level of background in linear algebra that 
readers interested in applications to medical imaging may not possess
so we feel that a more elementary presentation
might be a valuable supplement to  Higham \cite{HighamFM}. 
In this paper, I present a unified exposition of these results, 
including a proof of the existence of the  Real Jordan Form, 
and give more direct proofs of some of these
results using the Real Jordan Form.

\vfill\eject
\section{Introduction and Motivations}
\label{sec-intro}
Theorems about the conditions for the existence (and uniqueness)
of a real logarithm (or a real square root) of a real matrix
are the theoretical basis for 
various numerical methods for exponentiating a matrix or for computing
its logarithm using a method known as {\it scaling and squaring\/}
(resp.  {\it inverse scaling and squaring\/}).
Such methods play an important role in the
{\it log-Euclidean framework\/} 
due to Arsigny, Fillard, Pennec and Ayache
and its applications to medical imaging
\cite{Arsignythesis,ArsignyLogE1,ArsignySIAM,ArsignyLogE2}.

\medskip
The registration of medical images is an important and difficult
problem. The work described in Arsigny, Commowick, Pennec and Ayache 
\cite{ArsignyLEPT} (and Arsigny's thesis \cite{Arsignythesis}) 
makes an orginal and valuable contribution
to this problem by describing a method for parametrizing a class
of non-rigid deformations with a small number of degrees of freedom.
After a global affine alignment, this sort of parametrization
allows a finer local registration with very smooth transformations.
This type of parametrization is particularly well adpated to the 
registration of histological slices, see Arsigny, Pennec and Ayache
\cite{ArsignyLogE2}.

\medskip
The goal is to fuse some affine or rigid transformations in such a way that
the resulting transformation is invertible and smooth.
The direct approach which consists in blending
$N$ global affine or rigid transformations, $T_1, \ldots, T_N$ 
using weights, $w_1, \ldots, w_N$,  does not work
because the resulting transformation,
\[
T = \sum_{i = 1}^N w_i T_i,
\]
is not necessarily invertible. 
The purpose of the weights is to define the domain of influence
in space of each $T_i$. 

\medskip
The novel key idea is to 
associate to each rigid (or affine) transformation, $T$, of $\reals^n$,
a vector field, $V$, and to view $T$ as the diffeomorphism,
$\Phi^V_1$, corresponding to the time $t = 1$, where $\Phi^V_t$ is the
global flow associated with $V$. In other words, $T$ is
the result of integrating an ODE
\[
X' = V(X, t),
\]
starting with some initial condition, $X_0$, and $T = X(1)$.

\medskip
Now, it would be highly desirable if the vector field, $V$, did
not depend on the time parameter, and this is indeed possible for 
a large class of affine transformations, which is one of
the nice contributions of the work of 
Arsigny, Commowick, Pennec and Ayache 
\cite{ArsignyLEPT}.

\medskip
Recall that an affine transformation, $X \mapsto LX + v$,
(where $L$ is an $n\times n$ matrix and $X, v \in \reals^n$) 
can be conveniently represented as a linear transformation
from $\reals^{n+1}$ to itself if we write
\[
\binom{X}{1} \mapsto \amsmata{L}{v}{0}{1}\binom{X}{1}.
\]
Then, the ODE with constant coefficients
\[
X' = LX + v,
\]
can be written 
\[
\binom{X'}{0} = \amsmata{L}{v}{0}{0}\binom{X}{1}
\]
and, for every initial condition, $X = X_0$, its unique solution is given by
\[
\binom{X(t)}{1} = \exp\left(t\amsmata{L}{v}{0}{0} \right)\binom{X_0}{1}.
\]
Therefore, if we can find reasonable conditions on matrices,
$T =  \amsmata{M}{t}{0}{1}$, to ensure that they have  a 
unique real logarithm,
\[
\log(T) = \amsmata{L}{v}{0}{0},
\]
then we will be able to associate a vector field, $V(X) = LX + v$, to $T$,
in such a way that $T$ is recovered by integrating the ODE,
$X' = LX + v$. Furthermore, given $N$ transformations,
$T_1, \ldots, T_N$, such that $\log(T_1), \ldots, \log(T_N)$
are uniquely defined, we can fuse $T_1, \ldots, T_N$  at the
{\it infinitesimal level\/} by defining the ODE obtained by blending
the  vector fields, $V_1, \ldots, V_N$, associated with
$T_1, \ldots, T_N$ (with $V_i(X) = L_iX + v_i$), namely
\[
V(X) = \sum_{i = 1}^N w_i(X)(L_i X + v_i).
\]
Then, it is easy to see that the ODE, 
\[
X' = V(X),
\]
has a unique solution for every
$X = X_0$ defined for all $t$, and the fused transformation
is just $T = X(1)$. Thus, the fused vector field, 
\[
V(X)  = \sum_{i = 1}^N w_i(X)(L_i X + v_i),
\]
yields a one-parameter group of diffeomorphisms, $\Phi_t$.
Each transformation, $\Phi_t$, is smooth and invertible and
is called a {\it Log-Euclidean polyaffine tranformation\/}, for short,
{\it LEPT\/}.
Of course, we have the equation
\[
\Phi_{s + t} = \Phi_s \circ \Phi_t,
\]
for all $s, t\in \reals$ so, in particular, the inverse of 
$\Phi_t$ is $\Phi_{-t}$. We can also interpret $\Phi_s$ as
$(\Phi_1)^s$, which will yield a fast method for computing 
$\Phi_s$.
Observe that when the weight are scalars, the one-parameter group is given by
\[
\binom{\Phi_t(X)}{1} = 
\exp\left(t  \sum_{i = 1}^N w_i\amsmata{L_i}{v_i}{0}{0} \right)
\binom{X}{1},
\]
which is the Log-Euclidean mean of the affine transformations,
$T_i$'s (w.r.t. the weights $w_i$).

\medskip
Fortunately, there is a sufficient condition for a real matrix to have
a unique real logarithm and this condition is not too
restrictive in practice.

\medskip
Let  $\s{S}(n)$ denotes the set of all  real matrices
whose  eigenvalues, $\lambda + i\mu$, lie in the horizontal strip
determined by the condition
$-\pi < \mu < \pi$. We have the following weaker version of
Theorem \ref{logexistuniq}:

\begin{thm}
\label{logthm1}
The image, $\exp(\s{S}(n))$, of $\s{S}(n)$ by the exponential map
is the set of real invertible matrices with no negative eigenvalues
and $\mapdef{\exp}{\s{S}(n)}{\exp(\s{S}(n))}$ is a bijection. 
\end{thm}

\medskip
Theorem \ref{logthm1} is stated in Kenney and Laub  \cite{KenneyLaub}
without proof. Instead, Kenney and Laub cite 
DePrima and Johnson \cite{DePrimaJohnson} for a proof but this latter
paper deals with  complex matrices and does not contain a proof 
of our result either.

\medskip
It is also known that under the same condition
(no negative eigenvalues)
every real $n\times n$ matrix, $A$, has
a real square root, that is, there is a real matrix, $X$, such that
$X^2 = A$. Moreover, if the eigenvalues, $\rho\, e^{i\theta}$, of $X$
satisfy the condition $-\frac{\pi}{2} < \theta < \frac{\pi}{2}$, 
then $X$ is unique (see Theorem \ref{squareexistuniq}).

\medskip
Actually, there is a necessary and sufficient 
condition for a real matrix to have a real logarithm
(or a real square root) but it is fairly subtle
as it involves the parity of the number of Jordan blocks
associated with negative eigenvalues, see Theorem \ref{log2}.
The first occurrence of this theorem  that we have 
found in the literature is a paper by Culver \cite{Culver} published in
1966. We offer a proof using Theorem
\ref{Jordanform3}, which is more explicit than Culver's proof.

\medskip
Curiously, complete and direct proofs of the main Theorems,
\ref{log2}, \ref{logexistuniq}, and \ref{squareexistuniq}
do not seem to exist and references 
found in various papers are sometimes incorrect (for more on this,
see the beginning of Section \ref{sec2}, the remark after
the proof of Theorem \ref{sqmat2} and the remark after
the proof of Theorem \ref{squareexistuniq}).
Versions of these results do appear
in Higham's book \cite{HighamFM} but one of the theorems
involved (Theorem 1.28) is not proved and 
closer examination reveals that Theorem 1.36  (in Higham's book)
is needed to prove Theorem 1.28.

\medskip
In view of all this, we feel that providing a unifying treatment
and giving complete proofs of these results will be of value
to the mathematical community.

\section{Jordan Decomposition and the Real Jordan Form}
\label{sec1}
The proofs of the results stated in Section \ref{sec-intro}
make heavy use
of the Jordan normal form of a matrix and its cousin, the Jordan
decomposition of a matrix into its semisimple part and its nilpotent part.
The purpose of this section is to review these concepts rather
thoroughly to make sure that the reader has the background necessary
to understand the proofs in Section \ref{sec2} and  Section \ref{sec3}.
We pay particular attention to the {\it Real Jordan Form\/}
(Horn and Johnson \cite{HornJohn}, Chapter 3, Section 4, Theorem 3.4.5,
Hirsh and Smale \cite{HirshSmale} Chapter 6)
which, although familiar to experts in linear algebra, is typically missing 
from ``standard'' algebra books. We give a complete proof of
the Real Jordan Form as such a proof does not seem to be easily found
(even  Horn and Johnson \cite{HornJohn} only give a sketch of the proof,
but it is covered in Hirsh and Smale \cite{HirshSmale}, Chapter 6).

\medskip
Let $V$ be a finite dimensional real vector space.
Recall that we can form the {\it complexification\/}, $V_{\complex}$, of $V$.
The space  $V_{\complex}$ is the complex vector space, $V\times V$, with
the addition operation given by 
\[
(u_1, v_1) + (u_2, v_2) = (u_1 + u_2, v_1 + v_2), 
\]
and  the scalar multiplication given by
\[
(\lambda + i\mu)\cdot (u, v) = (\lambda u - \mu v, \mu u + \lambda v)
\qquad (\lambda, \mu\in \reals).
\]
Obviously
\[
(0, v) = i\cdot (v, 0), 
\]
so every vector, $(u, v)\in V_{\complex}$, can written uniquely as
\[
(u, v) = (u, 0) + i\cdot (v, 0).
\] 
The map  from $V$ to $V_{\complex}$ given by $u \mapsto (u, 0)$
is obviously an injection and for notational convenience,
we write $(u, 0)$ as $u$, we suppress the symbol (``dot'')
for scalar multiplication and we write
\[
(u, v) = u + i v, \qquad\hbox{with}\> u, v\in V.
\] 
Observe that if $(e_1, \ldots, e_n)$ is a basis of $V$, then
it is also a basis of $V_{\complex}$.

\medskip
Every linear map, $\mapdef{f}{V}{V}$, yields a linear map,
$\mapdef{f_{\complex}}{V_{\complex}}{V_{\complex}}$,  with
\[
f_{\complex}(u + i v ) = f(u) + i f(v), \qquad\hbox{for all}\> u, v\in V.
\]

\begin{defin}
\label{nilpotdef}
A linear map, $\mapdef{f}{V}{V}$, is {\it semisimple\/}
iff $f_{\complex}$ can be diagonalized.
In terms of matrices, a real matrix, $A$, is {\it semisimple\/} iff there 
are some matrices $D$ and $P$ with entries in $\complex$,
with $P$ invertible and $D$ a diagonal matrix, so that
$A = PD P^{-1}$.
We say that $f$ is {\it nilpotent\/} iff $f^r = 0$ for some
positive integer, $r$, and a matrix, $A$, is {\it nilpotent\/} iff
$A^r = 0$ for some  positive integer, $r$.
We say that $f$ is {\it unipotent\/} iff $f - \id$ is nilpotent
and a matrix $A$ is {\it unipotent\/} iff $A - I$ is nilpotent.
\end{defin}

\medskip
If $A$ is unipotent, then $A = I + N$ where $N$ is nilpotent.
If $r$ is the smallest integer so that $N^r = 0$ (the 
{\it index of nilpotency\/} of $N$),
then it is easy to check that
\[
I - N + N^2 + \cdots + (-1)^{r - 1} N^{r - 1}
\]
is the inverse of  $A = I + N$.

\medskip
For example, rotation matrices are semisimple, although in general
they can't be diagonalized over $\reals$, since their
eigenvalues are complex numbers of the form $e^{i\theta}$.
Every upper-triangular matrix where all the diagonal entries are  zero
is nilpotent.

\begin{defin}
\label{Jordandecomp}
If $\mapdef{f}{V}{V}$ is a linear map with $V$ a finite
vector space over $\reals$ or $\complex$, a {\it Jordan decomposition\/}
of $f$ is a pair of linear maps, $\mapdef{f_S, f_N}{V}{V}$,
with $f_S$ semisimple  and $f_N$ nilpotent, such that
\[
f = f_S + f_N\qquad
\hbox{and} \qquad
f_S\circ f_N = f_N\circ f_S.
\]
\end{defin}

\medskip
The theorem below is a very useful technical tool
for dealing with the exponential map.
It can be proved from the so-called primary decomposition theorem or 
from the Jordan form (see Hoffman and Kunze \cite{HoffmanKunze},
Chapter 6, Section 4 or
Bourbaki \cite{BourbakiA2}, Chapter VII, \S5).

\begin{thm}
\label{Jordform1}
If $V$ is a finite dimensional vector space over $\complex$, then
every linear map, $\mapdef{f}{V}{V}$, has a unique Jordan
decomposition, $f = f_S + f_N$. Furthermore, $f_S$ and $f_N$ can
be expressed as polynomials in $f$ with no constant term.
\end{thm}

\remark
In fact, Theorem \ref{Jordform1} holds for 
any finite dimensional vector space over a
{\it perfect field\/}, $K$ (this means that 
either $K$ has characteristic zero of that
$K^p = K$,
where $K^p = \{a^p \mid a\in K\}$ and where $p\geq 2$ is the
characteristic of the field $K$). 
The proof of this stronger version of Theorem \ref{Jordform1} 
is more subtle and involves some elementary Galois theory
(see Hoffman and Kunze \cite{HoffmanKunze},
Chapter 7, Section 4 or, for maximum generality,
Bourbaki \cite{BourbakiA2}, Chapter VII, \S5).

\medskip
We will need Theorem \ref{Jordform1} in the case where
$V$ is a real vector space. In fact we need a slightly
refined version of  Theorem \ref{Jordform1} for $K = \reals$
known as the {\it Real Jordan form\/}.
First, let us review  Jordan matrices and real Jordan matrices.

\begin{defin}
\label{Jordanformdef}
A (complex) {\it Jordan block\/} is an  $r \times r$ matrix, $J_r(\lambda)$,
of the form
\[
J_r(\lambda) = 
\begin{pmatrix}
\lambda & 1 & 0 &  \cdots &  0 \\
0 & \lambda  & 1  & \cdots & 0 \\
\vdots & \vdots  & \ddots  & \ddots  &  \vdots      \\
0 & 0 & 0 & \ddots &  1 \\
0 & 0 & 0 & \cdots &  \lambda
\end{pmatrix},
\]
where $\lambda\in \complex$, with $J_1(\lambda) = (\lambda)$ if $r = 1$.
A {\it real Jordan block\/} is either 
\begin{enumerate}
\item[(1)]
a Jordan block as above with $\lambda\in \reals$, or
\item[(2)] 
a real $2r \times 2r$ matrix, $J_{2r}(\lambda, \mu)$, of the form
\[
J_{2r}(\lambda, \mu) = 
\begin{pmatrix}
L(\lambda, \mu) & I & 0 &  \cdots &  0 \\
0 & L(\lambda, \mu)  & I  & \cdots & 0 \\
\vdots & \vdots  & \ddots  & \ddots  &  \vdots      \\
0 & 0 & 0 & \ddots &  I \\
0 & 0 & 0 & \cdots &  L(\lambda, \mu)
\end{pmatrix},
\]
where $L(\lambda, \mu)$ is a  $2\times 2$ matrix of the form
\[
L(\lambda, \mu) = 
\begin{pmatrix}
\lambda &  -\mu \\
\mu & \lambda
\end{pmatrix},
\]
with $\lambda, \mu\in \reals$, $\mu\not= 0$,  with
$I$ the $2\times 2$ identity matrix and with
$J_{2}(\lambda, \mu) = L(\lambda, \mu)$ when $r = 1$.
\end{enumerate}
A (complex) {\it Jordan matrix\/}, $J$, is an $n \times n$
block diagonal matrix of the form
\[
J = 
\begin{pmatrix}
J_{r_1}(\lambda_1) &  \cdots & 0 \\
 \vdots  & \ddots   & \vdots    \\
 0 & \cdots & J_{r_m}(\lambda_m)
\end{pmatrix},
\]
where each $J_{r_k}(\lambda_k)$ is a (complex) Jordan block 
associated with some $\lambda_k \in \complex$ and with \\
$r_1 + \cdots + r_m = n$.
A {\it real Jordan matrix\/}, $J$, is an $n \times n$
block diagonal matrix of the form
\[
J = 
\begin{pmatrix}
J_{s_1}(\alpha_1) &  \cdots & 0 \\
 \vdots  & \ddots   & \vdots    \\
 0 & \cdots & J_{s_m}(\alpha_m)
\end{pmatrix},
\]
where each $J_{s_k}(\alpha_k)$ is a real Jordan block 
either associated with some $\alpha_k = \lambda_k \in \reals$ as in (1)
or associated with some $\alpha_k = (\lambda_k, \mu_k) \in \reals^2$,
with $\mu_k \not= 0$, as in (2), in which case $s_k = 2r_k$.
\end{defin}

\medskip
To simplify notation, we often write $J(\lambda)$ for $J_r(\lambda)$
(or  $J(\alpha)$ for $J_s(\alpha)$).
Here is an example of a Jordan matrix with four blocks:
\[
J =
\begin{pmatrix}
\lambda & 1       & 0      &  0      & 0       &  0      & 0   & 0 \\
0       & \lambda & 1      &  0      & 0       &  0      & 0   & 0 \\
0       & 0       &\lambda &  0      & 0       &  0      & 0   & 0 \\
0       & 0       &  0     & \lambda & 1       &  0      & 0   & 0 \\
0       & 0       &  0     & 0       & \lambda &  0      & 0   & 0 \\
0       & 0       &  0     & 0       & 0       & \lambda & 0   & 0 \\
0       & 0       &  0     & 0       & 0       & 0       & \mu & 1 \\  
0       & 0       &  0     & 0       & 0       & 0       & 0   & \mu    
\end{pmatrix}.
\]

\medskip
In order to prove properties of the exponential of
Jordan blocks, we need to understand the deeper reasons
for the existence of the Jordan form. For this, we review the notion of
a minimal polynomial.

\medskip
Recall that a polynomial, $p(X)$, of degree $n \geq 1$  is a 
{\it monic polynomial\/} 
iff the monomial of highest degree in $p(X)$ is of the form $X^n$
(that is, the coefficient of $X^n$ is equal to $1$).
As usual, let $\complex[X]$ be the ring of polynomials
\[
p(X) = a_0X^n + a_1 X^{n-1} + \cdots + a_{n -1} X + a_n,
\]
with complex coefficient, $a_i \in \reals$, and let
$\reals[X]$ be the ring of polynomials with real coefficients, $a_i\in\reals$.
If $V$ is a finite dimensional complex vector space and $\mapdef{f}{V}{V}$
is a given linear map, every polynomial
\[
p(X) = a_0X^n + a_1 X^{n-1} + \cdots + a_{n -1} X + a_n,
\]
yields the linear map denoted $p(f)$, where
\[
p(f)(v) = a_0f^n(v) + a_1 f^{n-1}(v) + \cdots + a_{n -1} f(v) + a_nv,
\qquad\hbox{for every $v\in V$},
\] 
and where $f^k = f\circ \cdots \circ f$ is
the composition of $f$ with itself $k$ times.
We also write
\[
p(f) = a_0f^n + a_1 f^{n-1} + \cdots + a_{n -1} f + a_n\id.
\]

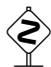
\danger
Do not confuse $p(X)$  and $p(f)$.
The expression $p(X)$ denotes a polynomial in the ``indeterminate'' $X$,
whereas $p(f)$ denotes a linear map from $V$ to $V$.

\medskip
For example, if $p(X)$ is the polynomial
\[
p(X) = X^3 - 2X^2 + 3X -1,
\]
if $A$ is any $n\times n$ matrix, then $p(A)$ is the $n\times n$ matrix
\[
p(A) = A^3 - 2A^2 + 3A -I
\]  
obtained by formally substituting the matrix $A$ for the variable $X$.

\medskip
Thus, we can define a ``scalar multiplication'',
$\mapdef{\cdot}{\complex[X]\times V}{V}$, by
\[
p(X) \cdot v = p(f)(v), \qquad v\in V.
\]
We immediately check that
\begin{eqnarray*}
p(X)\cdot (u + v) & = & p(X)\cdot u + p(X)\cdot v \\
(p(X) + q(X))\cdot u & = & p(X)\cdot u +  q(X)\cdot u \\  
(p(X)q(X))\cdot u & = & p(X)\cdot (q(X)\cdot u) \\
1\cdot u & = & u, 
\end{eqnarray*}
for all $u, v\in V$ and all $p(X), q(X)\in \complex[X]$,  
where $1$ denotes the polynomial of degree $0$ with constant term $1$.

\medskip
It follows that the scalar multiplication,
$\mapdef{\cdot}{\complex[X]\times V}{V}$, makes $V$ into a $\complex[X]$-module
that we will denote by $V_f$. Furthermore, as $\complex$ is a subring of
$\complex[X]$ and as $V$ is finite-dimensional, $V$ is finitely generated over
$\complex$ and so $V_f$ is finitely generated as a module over $\complex[X]$.

\medskip
Now, because $V$ is finite dimensional, we claim that there is some
polynomial, $q(X)$, that {\it annihilates\/} $V_f$, that is, so that
\[
q(f)(v) = 0, \qquad \hbox{for all} \quad v\in V.
\]
To prove this fact, observe that if $V$ has dimension $n$, then
the set of linear maps from $V$ to $V$ has dimension $n^2$.
Therefore any $n^2 + 1$ linear maps must be linearly dependent, so
\[
\id, f, f^2, \ldots, f^{n^2}
\]
are linearly dependent linear maps and
there is a nonzero polynomial, $q(X)$, of degree at most $n^2$ so that
$q(f)(v) = 0$ for all $v\in V$.
(In fact, by the {\it Cayley-Hamilton Theorem\/}, the characteristic polynomial,
$q_f(X) = \det(X\,\id  - f)$,  of $f$ annihilates $V_f$, 
so there is some annihilating polynomial
of degree at most $n$.)
By abuse of language (and notation),
if $q(X)$ annihilates $V_f$, we also say that $q(X)$
annihilates $V$.

\medskip
The set of annihilating polynomials of $V$ forms a principlal ideal
in $\complex[X]$, which
means that there is a unique monic polynomial of minimal degree, $p_f$,
annihilating $V$ and every other polynomial annihilating $V$ is a multiple of $p_f$.
We call this minimal monic polynomial  annihilating $V$ the 
{\it minimal polynomial\/} of $f$. 

\medskip
The fact that $V$ is annihilated by some polynomial in $\complex[X]$
makes $V_f$ a {\it torsion\/} $\complex[X]$-module. Furthermore,
the ring  $\complex[X]$ has the property that every ideal is a {\it principal
ideal domain\/}, abbreviated PID
(this means that every ideal is generated by a single
polynomial which can be chosen to be monic and of smallest degree).
The ring $\reals[X]$ is also a PID. 
In fact, the ring $k[X]$ is a PID for any field, $k$.
But then, we can apply some powerful results about the structure of
finitely generated torsion modules over a PID to $V_f$
and obtain various decompositions of $V$ into subspaces
which yield useful
normal forms for $f$, in particular, the Jordan form.

\medskip
Let us give one more definition
before stating our next  important theorem:
Say that $V$ is a {\it cyclic module\/} iff $V$ is generated by
a single element as a $\complex[X]$-module, which means that
there is some $u\in V$ so that
$u, f(u), f^2(u), \ldots, f^k(u), \ldots,$
generate $V$. 

\begin{thm}
\label{elemdiv1}
let $V$ be a finite-dimensional complex vector space of dimension $n$.
For every linear map, $\mapdef{f}{V}{V}$, there is a direct sum decomposition,
\[
V = V_1 \oplus V_2 \oplus \cdots \oplus V_m,
\]
where each $V_i$ is a cyclic $\complex[X]$-module 
such that the minimal polynomial
of the restriction of $f$ to $V_i$ is of the form $(X - \lambda_i)^{r_i}$.
Furthermore, the number, $m$, of subspaces $V_i$ 
and the minimal polynomials of the $V_i$ are uniquely determined by $f$ 
and, for each such polynomial, $(X - \lambda)^r$, the number, $m_i$, of
$V_i$'s that have  $(X - \lambda)^r$ as minimal polynomial
(that is, if $\lambda = \lambda_i$ and $r = r_i$) is  uniquely determined by $f$. 
\end{thm}

\medskip
A proof of Theorem \ref{elemdiv1} can be found in
M. Artin \cite{Artin91}, Chapter 12, Section 7,
Lang \cite{Lang93}, Chapter XIV, Section 2, 
Dummit and Foote \cite{DummitFoote}, Chapter 12, Section 1 and Section 3,
or D. Serre \cite{SerreDenis}, Chapter 6, Section 3.
A very good exposition is also given in
Gantmacher \cite{Gantmacher1}, Chapter VII, in particular, see Theorem 8
and Theorem 12. However, in Gantmacher,
elementary divisors are defined
in a rather cumbersome manner in terms of ratios of determinants 
of certain minors. This makes, at times, the proof unnecessarily hard 
to follow.

\medskip
The minimal polynomials, $(X - \lambda_i)^{r_i}$, associated
with the $V_i$'s are called the {\it elementary divisors\/} of $f$.
They need not be distinct. To be more precise, if
the set of distinct elementary divisors of $f$ is
\[
\{(X - \lambda_1)^{r_1}, \ldots, (X - \lambda_t)^{r_t}\} 
\]
then $(X - \lambda_1)^{r_1}$ appears $m_1 \geq 1$ times,
$(X - \lambda_2)^{r_2}$ appears $m_2 \geq 1$ times, ...,
$(X - \lambda_t)^{r_t}$ appears $m_t \geq 1$ times, with
\[
m_1 + m_2 + \cdots + m_t = m.
\]
The number, $m_i$, is called the {\it multiplicity\/} of $(X - \lambda_i)^{r_i}$.
Furthermore, if $(X - \lambda_i)^{r_i}$ and $(X - \lambda_j)^{r_j}$
are two distinct elementary divisors,  it is possible 
that $r_i \not= r_j$  yet  $\lambda_i = \lambda_j$.

\medskip
Observe that $(f - \lambda_i\id)^{r_i}$ is nilpotent on $V_i$
with index of nilpotency $r_i$ (which means that 
$(f - \lambda_i\id)^{r_i} = 0$ on $V_i$ but 
$(f - \lambda_i\id)^{r_i - 1}\not= 0$ on $V_i$). Also,
note that the monomials, $(X - \lambda_i)$, are the irreducible factors
of the minimal polynomial of $f$. 

\medskip
Next, let us take a closer look at the subspaces, $V_i$.
It turns out that we can find a ``good'' basis of $V_i$ 
so that in this basis, the restriction of $f$ to $V_i$ 
is a Jordan block.

\begin{prop}
\label{elemdiv2}
Let $V$ be a finite-dimensional vector space 
and let $\mapdef{f}{V}{V}$ be a linear map.
If $V$ is a cyclic $\complex[X]$-module and if $(X - \lambda)^n$ is
the minimal polynomial of $f$, then there is a basis of $V$ of the form
\[
((f - \lambda \id)^{n-1}(u), (f - \lambda \id)^{n-2}(u), \ldots,
(f - \lambda \id)(u), u),
\]
for some $u\in V$. With respect to this basis, the matrix of $f$
is the Jordan block
\[
J_n(\lambda) = 
\begin{pmatrix}
\lambda & 1 & 0 &  \cdots &  0 \\
0 & \lambda  & 1  & \cdots & 0 \\
\vdots & \vdots  & \ddots  & \ddots  &  \vdots      \\
0 & 0 & 0 & \ddots &  1 \\
0 & 0 & 0 & \cdots &  \lambda
\end{pmatrix}.
\]
Consequently, $\lambda$ is an eigenvalue of $f$.
\end{prop}

\begin{proof}
A proof is given in Section \ref{sec6}.
\end{proof}

\medskip
Using Theorem \ref{elemdiv1}  and Proposition \ref{elemdiv2}
we get the Jordan form for complex matrices.

\begin{thm} (Jordan Form)
\label{Jordanform1}
For every complex $n\times n$ matrix, $A$, there is some
invertible matrix, $P$, and some Jordan matrix, $J$, so that
\[
A = P J P^{-1}.
\]
If $\{\lambda_1, \ldots, \lambda_s\}$ is the set of 
eigenvalues of $A$, then the diagonal elements of the
Jordan blocks of $J$ are among the $\lambda_i$ and
every $\lambda_i$ corresponds to
one of more Jordan blocks of $J$.
Furthermore, the number, $m$, of Jordan blocks, 
the distinct Jordan block, $J_{r_{i}}(\lambda_i)$,  and
the number of times, $m_i$, that each Jordan block, $J_{r_{i}}(\lambda_i)$, occurs
are uniquely determined by $A$.
\end{thm}

\medskip
The number $m_i$ is called the {\it multiplicity\/} of the block 
$J_{r_{i}}(\lambda_i)$. Observe that the column vector associated with
the first entry of every Jordan block is an eigenvector of $A$.
Thus, the number, $m$, of Jordan blocks is the number of linearly
independent eigenvectors of $A$.

\medskip
Beside the references that we cited for the proof
of Theorem \ref{elemdiv1}, other proofs
of Theorem \ref{Jordanform1} can be found in the literature.
Often, these proofs do not cover the uniqueness statement.
For example, a nice proof  is given in Godement 
\cite{Godementalg}, Chapter 35.
Another interesting proof is given in Strang \cite{Strang88}, Appendix B.
A more  ``computational proof'' is given in
Horn and Johnson, \cite{HornJohn}, Chapter 3, Sections 1-4.

\medskip
Observe that Theorem \ref{Jordanform1} implies that
the characteristic polynomial, $q_f(X)$, of $f$ is the product
of the elementary divisors of $f$ (counted with their multiplicity).
But then, $q_f(X)$ must annihilate $V$. Therefore, we obtain
a quick proof of the Cayley Hamilton Theorem (of course, we had to work
hard to get Theorem \ref{Jordanform1}!). 
Also, the minimal polynomial of $f$ is the 
least common multiple (lcm) of the elementary divisors of $f$.

\medskip
The following technical result will be needed for finding
the logarithm of a real matrix:

\begin{prop}
\label{rearrange}
If $J$ is a $2n \times 2n$ complex Jordan matrix consisting
of two conjugate blocks $J_n(\lambda + i\mu)$ and
$J_n(\lambda - i\mu)$ of dimension $n$ ($\mu\not= 0$), then
there is a permutation matrix, $P$, and matrix, $E$,  so that
\[
J = P E P^{-1},
\]
where $E$ is a block matrix of the form
\[
E = 
\begin{pmatrix}
D & I & 0 &  \cdots &  0 \\
0 & D  & I  & \cdots & 0 \\
\vdots & \vdots  & \ddots  & \ddots  &  \vdots      \\
0 & 0 & 0 & \ddots &  I \\
0 & 0 & 0 & \cdots &  D
\end{pmatrix},
\]
and with $D$  the diagonal $2\times 2$ matrix
\[
D = 
\begin{pmatrix}
\lambda + i\mu &  0 \\
0 & \lambda - i\mu
\end{pmatrix}.
\]
Furthermore, there is a complex invertible matrix, $Q$, and a real
Jordan matrix,  $C$,
so that
\[
J = Q C Q^{-1},
\]
where $C$ is of the form
\[
C = 
\begin{pmatrix}
L & I & 0 &  \cdots &  0 \\
0 & L  & I  & \cdots & 0 \\
\vdots & \vdots  & \ddots  & \ddots  &  \vdots      \\
0 & 0 & 0 & \ddots &  I \\
0 & 0 & 0 & \cdots &  L
\end{pmatrix},
\]
with 
\[
L = 
\begin{pmatrix}
\lambda &  -\mu \\
\mu & \lambda
\end{pmatrix}.
\]
\end{prop}

\begin{proof}
First, consider an example, namely,
\[
J = 
\begin{pmatrix}
\lambda + i\mu & 1 & 0 & 0 \\
0 & \lambda + i\mu & 0 & 0  \\
0 & 0 & \lambda - i\mu & 1  \\ 
0 & 0 & 0 & \lambda - i\mu 
\end{pmatrix}.
\]
If we permute rows $2$ and $3$, we get
\[
\begin{pmatrix}
\lambda + i\mu & 1 & 0 & 0 \\
0 & 0 & \lambda - i\mu & 1  \\ 
0 & \lambda + i\mu & 0 & 0  \\
0 & 0 & 0 & \lambda - i\mu 
\end{pmatrix}
\]
and we permute columns $2$ and $3$, we get our matrix,
\[
E = \begin{pmatrix}
\lambda + i\mu & 0 & 1 & 0 \\
0 &  \lambda - i\mu & 0 & 1 \\ 
0 &  0 & \lambda + i\mu &  0  \\
0 & 0 & 0 & \lambda - i\mu 
\end{pmatrix}.
\]
We leave it as an exercise to generalize this method to two $n\times n$ 
conjugate Jordan blocks
to prove that we can find a permutation matrix, $P$, so that
$E = P^{-1} J P$ and thus, $J = P E P^{-1}$.

\medskip
Next, as $\mu\not= 0$, the matrix $L$ can be diagonalized and one
easily checks that
\[
D =
 \begin{pmatrix}
\lambda +i\mu & 0 \\
0  & \lambda - i\mu
\end{pmatrix} = 
 \begin{pmatrix}
-i & 1  \\
-i & -1
\end{pmatrix} 
\begin{pmatrix}
\lambda &  -\mu \\
\mu & \lambda
\end{pmatrix}
 \begin{pmatrix}
-i & 1  \\
-i & -1
\end{pmatrix}^{-1}.
\]
Therefore, using the block diagonal matrix
$S = \mathrm{diag}(S_2, \ldots, S_2)$ 
consisting of $n$ blocks
\[
S_2 =  \begin{pmatrix}
-i & 1  \\
-i & -1
\end{pmatrix},  
\]
we see that
\[
E = S C S^{-1}
\]
and thus,
\[
J =  P S C  S^{-1} P^{-1},
\]
which yields our second result with $Q =  P S$.
\end{proof}

\medskip
Proposition \ref{rearrange} shows that  every (complex) matrix, $A$,
is similar to a real Jordan matrix. Unfortunately, if $A$ is a {\it real\/}
matrix, there is no guarantee that we can find a {\it real\/} invertible
matrix, $P$, so that $A = P J P^{-1}$, with $J$ a real Jordan matrix.
This result known as the {\it Real Jordan Form\/}
is actually true but requires some work to be established. 
In this section, we state the theorem without proof.
A proof based on Theorem \ref{elemdiv1} is given in Section \ref{sec6}.

\begin{thm} (Real Jordan Form)
\label{Jordanform2}
For every real $n\times n$ matrix, $A$, there is some
invertible (real) matrix, $P$, and some real Jordan matrix, $J$, so that
\[
A = P J P^{-1}.
\]
For every Jordan block, $J_r(\lambda)$, of type (1),
$\lambda$ is some real eigenvalue of $A$ and for
every Jordan block, $J_{2r}(\lambda, \mu)$, of type (2),
$\lambda + i \mu$ is a complex eigenvalue of $A$
(with $\mu\not= 0$). Every eigenvalue of $A$ corresponds to
one of more Jordan blocks of $J$.
Furthermore, the number, $m$, of Jordan blocks, 
the distinct Jordan block, $J_{s_i}(\alpha_i)$,  and
the number of times, $m_i$, that each Jordan block, $J_{s_i}(\alpha_i)$, occurs
are uniquely determined by $A$.
\end{thm}

\medskip
Let $A$ be a real matrix and let 
$(X - \alpha_1)^{r_1}, \ldots, (X - \alpha_m)^{m_1}$
be its list of elementary divisors or, equivalently,
let $J_{r_1}(\alpha_1), \ldots, J_{r_m}(\alpha_m)$
be its list of Jordan blocks.
If, for every $r_i$ and every real eigenvalue $\lambda_i < 0$, 
the number, $m_i$,  of Jordan blocks identical to $J_{r_i}(\alpha_i)$
is even, then there is a way to rearrange these blocks
using the technique of Proposition \ref{rearrange}
to obtain a version of the real Jordan form that makes it
easy to find logarithms (and square roots) of real matrices.

\begin{thm} (Real Jordan Form, Special Version)
\label{Jordanform3}
Let $A$ be a real $n\times n$ matrix and let 
$(X - \alpha_1)^{r_1}, \ldots, (X - \alpha_m)^{m_1}$
be its list of elementary divisors or, equivalently,
let $J_{r_1}(\alpha_1), \ldots$, $J_{r_m}(\alpha_m)$
be its list of Jordan blocks.
If, for every $r_i$ and every real eigenvalue $\alpha_i < 0$, 
the number, $m_i$,  of Jordan blocks identical to $J_{r_i}(\alpha_i)$
is even, then there is a real invertible matrix, $P$, and
a real Jordan matrix, $J'$, such that $A = PJ'P^{-1}$ and 
\begin{enumerate}
\item[(1)]
Every block, $J_{r_i}(\alpha_i)$, of $J$ 
for which $\alpha_i\in \reals$ and $\alpha_i\geq 0$ 
is a Jordan block of type (1) of $J'$
(as in Definition \ref{Jordanformdef}), or
\item[(2)]
For every block, $J_{r_i}(\alpha_i)$, of $J$
for which either $\alpha_i\in \reals$ and $\alpha_i < 0$
or $\alpha_i = \lambda_i + i\mu_i$ with $\mu_i\not= 0$
($\lambda_i, \mu_i\in \reals$), the corresponding
real Jordan block of $J'$ is defined as follows: 
\begin{enumerate}
\item[(a)]
If $\mu_i\not= 0$, then $J'$ contains
the real Jordan block $J_{2r_i}(\lambda_i, \mu_i)$
of type (2)  (as in Definition \ref{Jordanformdef}), or
\item[(b)]
If $\alpha_i < 0$ then $J'$ contains the real Jordan block
$J_{2r_i}(\alpha_i, 0)$ whose diagonal blocks are of the form
\[
L(\alpha_i, 0) = 
\begin{pmatrix}
\alpha_i &  0 \\
0  & \alpha_i
\end{pmatrix}.
\]
\end{enumerate}
\end{enumerate}
\end{thm}

\begin{proof}
By hypothesis, for every real eigenvalue, $\alpha_i < 0$,
for every $r_i$, the Jordan block, $J_{r_i}(\alpha_i)$, occurs an even number
of times say $2 t_i$, so by using a permutation, we may assume that
we have $t_i$ pairs of identical blocks $(J_{r_i}(\alpha_i), J_{r_i}(\alpha_i))$.
But then, for each pair of blocks of this form,
we can apply part (1) of 
Proposition \ref{rearrange} (since $\alpha_i$ is its own conjugate),
which yields our result.
\end{proof}

\medskip
\remark
The above result generalizes the fact that
when we have a rotation matrix, $R$,
the eigenvalues $-1$ occurring in the real block diagonal form of $R$
can be paired up.

\medskip
The following theorem shows that the ``structure'' of the Jordan form of
a  matrix is preserved under exponentiation. This is an
important result that will be needed
to establish the necessity of the criterion for a real matrix to have a real
logarithm. Again, in this section, we state the theorem without proof.
A proof is given in Section \ref{sec6}.

\begin{thm}
\label{expJord}
For any (real or complex) $n\times n$ matrix, $A$, if
$A = P J P^{-1}$ where $J$ is a Jordan matrix of the form
\[
J = 
\begin{pmatrix}
J_{r_1}(\lambda_1) &  \cdots & 0 \\
 \vdots  & \ddots   & \vdots    \\
 0 & \cdots & J_{r_m}(\lambda_m)
\end{pmatrix},
\]
then there is some invertible matrix, $Q$, so that
the  Jordan form of $e^A$ is given by
\[
e^A = Q\, e(J)\,  Q^{-1},
\]
where $e(J)$ is the Jordan matrix
\[
e(J) = 
\begin{pmatrix}
J_{r_1}(e^{\lambda_1}) &  \cdots & 0 \\
 \vdots  & \ddots   & \vdots    \\
 0 & \cdots & J_{r_m}(e^{\lambda_m})
\end{pmatrix},
\]
that is, each $J_{r_k}(e^{\lambda_k})$ is obtained from $J_{r_k}(\lambda_k)$
by replacing all the diagonal entries $\lambda_k$ by $e^{\lambda_k}$.
Equivalently, if the list of elementary divisors of $A$ is
\[
(X - \lambda_1)^{r_1}, \ldots, (X - \lambda_m)^{r_m},
\]
then the list of elementary divisors of $e^A$ is
\[
(X - e^{\lambda_1})^{r_1}, \ldots, (X - e^{\lambda_m})^{r_m}.
\]
\end{thm}

\section{Logarithms of Real Matrices; 
Criteria for Existence and Uniqueness}
\label{sec2}
If $A$ is any (complex) $n\times n$ matrix we say that a matrix, $X$, is
a {\it logarithm of $A$\/} iff $e^X = A$. 
Our goal is to find conditions for the existence and uniqueness
of real logarithms of real matrices.
The two main theorems of this section are Theorem \ref{log2}
and Theorem \ref{logexistuniq}. These theorems are used 
in  papers presenting methods for computing the 
logarithm of a matrix, including
Cheng,  Higham, Kenney and  Laub \cite{Higham01} and 
Kenney and Laub \cite{KenneyLaub}.

\medskip
Reference \cite{Higham01} 
cites  Kenney and Laub \cite{KenneyLaub} for a proof of Theorem 
\ref{logexistuniq}
but in fact, that paper does not give a proof.
Kenney and Laub \cite{KenneyLaub} do state  Theorem \ref{logexistuniq} 
as Lemma A2 of Appendix A, 
but they simply say that ``the proof is similar to that
of Lemma A1''. As to the proof of Lemma A1, Kenney and Laub
state without detail that it makes use of the Cauchy integral formula for
operators, a method used by
DePrima and Johnson \cite{DePrimaJohnson} to prove a similar theorem
for complex matrices (Section 4, Lemma 1) and where uniqueness
is also proved.
Kenney and Laub point out that the third hypothesis in that lemma is redundant.
Theorem \ref{logexistuniq} also appears in Higham's book \cite{HighamFM}
as Theorem 1.31. Its proof relies on Theorem 1.28 and
Theorem 1.18 (both in Higham's book)
but  Theorem 1.28 is not proved and 
only part of theorem 1.18 is proved in the text 
(closer examination reveals that Theorem 1.36  (in Higham's book)
is needed to prove Theorem 1.28).
Although Higham's Theorem 1.28 implies the injectivity statement
of  Theorem \ref{stardiffeo}
we feel that the proof of Theorem \ref{stardiffeo} is of
independent interest. Furthermore, Theorem \ref{stardiffeo} 
is a stronger result (it shows that $\exp$ is a diffeomorphism).

\medskip
Given this state of affairs where no explicit proof of 
Theorem \ref{logexistuniq} 
seems easily available,  we provide a complete proof of 
Theorem  \ref{logexistuniq}
using our special form of the Real Jordan Form.

\medskip
First, let us consider the case where $A$ is a complex matrix.
Now, we know that if $A = e^X$, then $\det(A) = e^{\mathrm{tr}(X)}\not= 0$, so
$A$ must be invertible. It turns out that
this condition is also sufficient.

\medskip
Recall that for every invertible matrix, $P$, and every matrix, $A$,
\[
e^{P A P^{-1}} = P e^A P^{-1} 
\]
and that for every block diagonal matrix, 
\[
A = 
\begin{pmatrix}
A_{1} &  \cdots & 0 \\
 \vdots  & \ddots   & \vdots    \\
 0 & \cdots & A_{m}
\end{pmatrix},
\]
we have
\[
e^A = 
\begin{pmatrix}
e^{A_{1}} &  \cdots & 0 \\
 \vdots  & \ddots   & \vdots    \\
 0 & \cdots & e^{A_{m}}
\end{pmatrix}.
\]
Consequenly, the problem of finding the logarithm of a matrix
reduces to the problem of finding the logarithm of a Jordan block $J_r(\alpha)$
with $\alpha\not= 0$.
However, every such Jordan block, $J_r(\alpha)$, can be written as
\[
J_r(\alpha) = \alpha I + H = \alpha I(I + \alpha^{-1} H),
\]
where $H$ is the nilpotent matrix
of index of nilpotency, $r$, given by  
\[
H = 
\begin{pmatrix}
0 & 1 & 0 &  \cdots &  0 \\
0 & 0  & 1  & \cdots & 0 \\
\vdots & \vdots  & \ddots  & \ddots  &  \vdots      \\
0 & 0 & 0 & \ddots &  1 \\
0 & 0 & 0 & \cdots &  0
\end{pmatrix}.
\] 
Furthermore, it is obvious that $N = \alpha^{-1} H$ is also nilpotent
of index of nilpotency, $r$, and we have
\[
J_r(\alpha) = \alpha I (I + N).
\]

\medskip
Logarithms of the diagonal matrix,
$\alpha I$, are easily found. If we write $\alpha = \rho e^{i\theta}$
where $\rho > 0$, then $\log \alpha = \log\rho + i (\theta +   2\pi h)$,
for any $h\in \integs$, and we can pick a logarithm of $\alpha I$ to be
\[
S = 
\begin{pmatrix}
\log \rho + i\theta & 0 & \cdots & 0 \\
0 & \log \rho + i\theta  &  \cdots & 0 \\
\vdots & \vdots  & \ddots   &  \vdots      \\
0 & 0 &  \cdots &  \log \rho + i\theta
\end{pmatrix}.
\]
Observe that if we can find a logarithm, $M$, of $I + N$, as $S$ commutes
with any matrix and as $e^{S} = \alpha I$ and $e^{M} = I + N$, we have
\[
e^{S + M} = e^{S}e^{M} = \alpha I(I + N) = J_r(\alpha),
\]
which means that $S + M$ is a logarithm of $J_r(\alpha)$.
Therefore, the problem reduces to finding the logarithm of a unipotent
matrix, $I + N$. However, this problem always has a solution.
To see this, remember that for $|u| < 1$, the power series
\[
\log(1 + u) = u - \frac{u^2}{2} + \frac{u^3}{3} + \cdots + 
 (-1)^{n +1} \frac{u^n}{n} + \cdots
\] 
is normally convergent. It turns out that the above fact can be generalized
to matrices in the following way:

\begin{prop}
\label{matrixlog}
For every $n\times n$ matrix, $A$, such that $\norme{A} < 1$, the series
\[
\log(I + A) = A - \frac{A^2}{2} + \frac{A^3}{3} + \cdots + 
 (-1)^{n +1} \frac{A^n}{n} + \cdots
\] 
is normally convergent for any matrix norm $\norme{\>}$ 
(a matrix norm satisfies the inequality
$\norme{AB} \leq \norme{A}\norme{B}$).
Furthermore, 
if $\norme{A} < 1$, then 
\[
e^{\log(I + A)} = I + A.
\]
\end{prop}

\medskip
\remark
For any matrix norm  $\norme{\>}$ and any
complex $n\times n$ matrix $A$, it can be shown that
\[
\rho(A) = \max_{1 \leq i \leq n} |\lambda_i| \leq \norme{A},
\]
where the $\lambda_i$ are the eigenvalues of $A$.
Furthermore, the set of (complex) diagonalizable matrices is dense
in the set of all complex matrices (see Serre \cite{SerreDenis}).
Using these two facts, it can be shown 
that if $\norme{A} < 1$, then 
\[
e^{\log(I + A)} = I + A
\]
for any matrix norm. 

\medskip
For any given $r \geq 1$,
the exponential and the logarithm (of matrices)
turn out to give a homeomorphim between the set of nilpotent matrices, $N$, 
and the set of unipotent matrices, $I + N$, for which $N^r = 0$. 
Let $\s{N}\mathit{il}(r)$ denote the
set of (real or complex)
nilpotent matrices of any dimension $n\geq 1$ such that
$N^r = 0$ and $\s{U}\mathit{ni}(r)$
denote the set of unipotent matrices, $U = I + N$, 
where $N \in  \s{N}\mathit{il}(r)$.
If $U = I + N\in  \s{U}\mathit{ni}(r)$,
note that $\log(I + N)$ is well-defined since
the power series for $\log(I + N)$ only 
has $r - 1$ nonzero terms,
\[
\log(I + N) = N - \frac{N^2}{2} + \frac{N^3}{3} + \cdots + 
 (-1)^{r} \frac{N^{r-1}}{r-1}.
\] 

\begin{prop}
\label{expnil}
The exponential map, $\mapdef{\exp}{\s{N}\mathit{il}(r)}{\s{U}\mathit{ni}(r)}$,
is a homeomorphism whose inverse is the logarithm.
\end{prop}

\begin{proof}
A complete proof can be found in Mmeimn\'e and Testard \cite{Mneimne},
Chapter 3, Theorem 3.3.3. The idea is to prove that
\[
\log(e^N) = N,
\>\hbox{for all}\> N\in \s{N}\mathit{il}(r)
\quad{and}\quad
e^{\log(U)} = U,
\>\hbox{for all}\> U\in \s{U}\mathit{ni}(r).
\]
To prove the first identity, it is enough to show that
for any fixed $N\in \s{N}\mathit{il}(r)$, we have
\[
\log(e^{t N}) = t N, 
\qquad \hbox{for all $t\in \reals$}.
\]
To do this, observe that the
functions $t \mapsto t N$ and $t \mapsto \log(e^{t N})$
are both equal to $0$ for $t = 0$. Thus, it is enough to show that their
derivatives are equal, which is left as an exercise.

\medskip
Next, for any $N\in \s{N}\mathit{il}(r)$,  the map
\[
t \mapsto e^{\log(I + t N)} - (I + tN),
\qquad t\in \reals
\]
is a polynomial, since $N^r = 0$. Furthermore, for $t$ sufficiently small,
$\norme{t N} < 1$ and in view of Proposition \ref{matrixlog},
we have $e^{\log(I + tN)} = I + t N$, so the above polynomial
vanishes in a neighborhood of $0$, which implies that it is 
identically zero. Therefore,
$e^{\log(I + N)} = I + N$, as required.
The continuity of $\exp$  and $\log$ is obvious.
\end{proof}

\medskip
Proposition \ref{expnil} shows that every unipotent matrix,
$I + N$, has the unique logarithm
\[
\log (I + N) = N - \frac{N^2}{2} + \frac{N^3}{3} + \cdots + 
 (-1)^{r} \frac{N^{r-1}}{r-1},
\]
where $r$ is the index of nilpotency of $N$. 
Therefore, if we let $M = \log(I + N)$, 
we have finally found a logarithm, $S + M$, for our original matrix, $A$.
As a result of all this, we have proved the following theorem:

\begin{thm}
\label{log1}
Every $n\times n$ invertible complex matrix, $A$, has a logarithm, $X$.
To find such a logarithm, we can proceed as follows:
\begin{enumerate}
\item[(1)]
Compute a Jordan form, $A = P J P^{-1}$, for $A$
and let $m$ be the number of Jordan blocks in $J$.
\item[(2)]
For every Jordan block, $J_{r_k}(\alpha_k)$, of $J$, write
$J_{r_k}(\alpha_j) = \alpha_k I(I + N_k)$, where $N_k$ is nilpotent.
\item[(3)]
If $\alpha_k = \rho_k e^{i\theta_k}$, with $\rho_k > 0$, let
\[
S_k = 
\begin{pmatrix}
\log \rho_k + i\theta_k & 0 & \cdots & 0 \\
0 & \log \rho_k + i\theta_k  &  \cdots & 0 \\
\vdots & \vdots  & \ddots   &  \vdots      \\
0 & 0 &  \cdots &  \log \rho_k + i\theta_k
\end{pmatrix}.
\]
We have $\alpha_k I = e^{S_k}$.
\item[(4)]
For every $N_k$, let
\[
M_k =  N_k - \frac{N_k^2}{2} + \frac{N_k^3}{3} + \cdots + 
 (-1)^{r_k} \frac{N^{r_k-1}}{r_k-1}.
\]
We have $I + N_k = e^{M_k}$.
\item[(5)]
If $Y_k = S_k + M_k$ and $Y$ is the block diagonal matrix
$\mathrm{diag}(Y_1, \ldots, Y_m)$, then
\[
X = P Y P^{-1}
\]
is a logarithm of $A$.
\end{enumerate}
\end{thm}

\medskip
Let us now assume that $A$ is a real matrix and let us try to find a
real logarithm. There is no problem in finding real logarithms of the
nilpotent parts but we run into trouble whenever an eigenvalue
is complex or real negative. Fortunately, we can circumvent these
problems by using the real Jordan form, provided that the condition
of Theorem \ref{Jordanform3} holds.

\medskip
The  theorem below gives a necessary and sufficient condition for a
real matrix to have a real logarithm. The first occurrence of this
theorem  that we have 
found in the literature is a paper by Culver \cite{Culver} published in
1966. The proofs in this paper rely heavily on results from
Gantmacher \cite{Gantmacher1}. Theorem \ref{log2} is also stated in
Horn and Johnson \cite{HornJohn2} as Theorem 6.4.15 (Chapter 6), but the proof
is left as an exercise. We offer a proof using Theorem
\ref{Jordanform3} which is more explicit than Culver's proof.

\begin{thm}
\label{log2}
Let $A$ be a real $n\times n$ matrix and let 
$(X - \alpha_1)^{r_1}, \ldots, (X - \alpha_m)^{m_1}$
be its list of elementary divisors or, equivalently,
let $J_{r_1}(\alpha_1), \ldots$, $J_{r_m}(\alpha_m)$
be its list of Jordan blocks. Then, $A$ has a real logarithm 
iff $A$ is invertible and if,
for every $r_i$ and every real eigenvalue $\alpha_i < 0$, 
the number, $m_i$,  of Jordan blocks identical to $J_{r_i}(\alpha_i)$
is even.
\end{thm}

\begin{proof}
First, assume that $A$ satisfies the conditions of Theorem \ref{log2}.
Since the matrix $A$ satisfies the condition
of  Theorem \ref{Jordanform3}, there is a real invertible matrix, $P$,
and a real Jordan matrix, $J'$, so that 
\[
A = P J' P^{-1},
\]
where $J'$ satisfies  conditions (1) and (2) of Theorem  \ref{Jordanform3}.
As $A$ is invertible, every block of $J'$ of the form $J_{r_k}(\alpha_k)$
corresponds to a real eigenvalue with  $\alpha_k > 0$ and we can write
$J_{r_k}(\alpha_j) = \alpha_k I(I + N_k)$, where $N_k$ is nilpotent.
As in Theorem \ref{log1} (4), we can find a real logarithm, $M_k$, of $I + N_k$
and as  $\alpha_k > 0$, the diagonal matrix  $\alpha_k I$
has the real logarithm
\[
S_k = 
\begin{pmatrix}
\log \alpha_k & 0 & \cdots & 0 \\
0 & \log \alpha_k   &  \cdots & 0 \\
\vdots & \vdots  & \ddots   &  \vdots      \\
0 & 0 &  \cdots &  \log \alpha_k 
\end{pmatrix}.
\]
Set $Y_k = S_k + M_k$.

\medskip
The other real Jordan blocks of $J'$ are of the form
$J_{2r_k}(\lambda_k, \mu_k)$,
with $\lambda_k, \mu_k \in \reals$, not both zero. Consequently,
we can write
\[
J_{2r_k}(\lambda_k, \mu_k) =  D_k + H_k = D_k(I + D_k^{-1}H_k) 
\]
where 
\[
D_k = 
\begin{pmatrix}
L(\lambda_k, \mu_k) &  \cdots & 0 \\
 \vdots  & \ddots   & \vdots    \\
 0 & \cdots & L(\lambda_k, \mu_k)
\end{pmatrix}
\]
with
\[
L(\lambda_k, \mu_k) = 
\begin{pmatrix}
\lambda_k & - \mu_k \\
\mu_k & \lambda_k
\end{pmatrix},
\]
and $H_k$ is a real nilpotent matrix. If we let $N_k = D_k^{-1}H_k$,
then $N_k$ is also nilpotent, $J_{2r_k}(\lambda_k, \mu_k) = D_k(I + N_k)$,
and we can find a logarithm, $M_k$, of $I + N_k$ as in Theorem \ref{log1} (4).
We can write $\lambda_k + i \mu_k = \rho_k e^{i\theta_k}$, with $\rho_k > 0$ and 
$\theta_k\in [-\pi, \pi)$, and then
\[
L(\lambda_k, \mu_k) = 
\begin{pmatrix}
\lambda_k & - \mu_k \\
\mu_k & \lambda_k
\end{pmatrix}
=
\rho_k
\begin{pmatrix}
\cos\theta_k & - \sin\theta_k \\
\sin\theta_k & \cos\theta_k
\end{pmatrix}.
\] 
If we set 
\[
S(\rho_k, \theta_k) = 
\begin{pmatrix}
\log\rho_k & - \theta_k \\
\theta_k & \log\rho_k
\end{pmatrix},
\]
a real matrix, we claim that
\[
L(\lambda_k, \mu_k) = e^{S(\rho_k, \theta_k)}.
\]
Indeed, $S(\rho_k, \theta_k) = \log\rho_k I + \theta_k E_2$,
with
\[
E_2 = 
\begin{pmatrix}
0 & -1 \\
1 & 0
\end{pmatrix},
\]
and it is well known that
\[
e^{\theta_k E_2} =
\begin{pmatrix}
\cos\theta_k & - \sin\theta_k \\
\sin\theta_k & \cos\theta_k
\end{pmatrix},
\]
so, as $\log\rho_k I$ and $\theta_k E_2$ commute, we get
\[
e^{S(\rho_k, \theta_k)} = e^{\log\rho_k I + \theta_k E_2} =
e^{\log\rho_k I}  e^{\theta_k E_2} = 
\rho_k \begin{pmatrix}
\cos\theta_k & - \sin\theta_k \\
\sin\theta_k & \cos\theta_k
\end{pmatrix}
= L(\lambda_k, \mu_k).
\]
If we form the real block diagonal matrix,
\[
S_k = 
\begin{pmatrix}
S(\rho_k, \theta_k) &  \cdots & 0 \\
 \vdots  & \ddots   & \vdots    \\
 0 & \cdots & S(\rho_k, \theta_k)
\end{pmatrix},
\]
we have $D_k = e^{S_k}$. Since $S_k$ and $M_k$ commute 
(observe that $M_k$ is obtained from adding up powers
of $N_k$ and $N_k$ only has $2\times 2$ blocks above a diagonal
of $2\times 2$ blocks and so, it commutes with a block diagonal matrix
of  $2\times 2$ blocks) and
\[
e^{S_k + M_k} = e^{S_k} e^{M_k} = D_k(I + N_k) = J_{2r_k}(\lambda_k, \mu_k),
\]
the matrix $Y_k = S_k + M_k$ is a logarithm
of $J_{2r_k}(\lambda_k, \mu_k)$.
Finally, if $Y$ is the block diagonal matrix
$\mathrm{diag}(Y_1, \ldots, Y_m)$, then
$X = P Y P^{-1}$ is a logarithm of $A$.

\medskip
Let us now prove that if $A$ has a real logarithm, $X$, then $A$ satisfies
the condition of Theorem \ref{log2}. 
As we said before, $A$ must be invertible. Since $X$ is a real matrix,
we know from the proof of Theorem \ref{Jordanform2} that
the Jordan blocks of $X$ associated with complex eigenvalues occur
in conjugate pairs, so they are of the form
\begin{align*}
& J_{r_k}(\alpha_k), \quad \alpha_k \in \reals, \\
& J_{r_k}(\alpha_k)\quad\hbox{and}\quad  J_{r_k}(\overline{\alpha}_k),
\quad \alpha_k = \lambda_k + i\mu_k,\> \mu_k \not= 0.
\end{align*}
By Theorem \ref{expJord}, the Jordan blocks of $A = e^X$ are 
obtained by replacing each $\alpha_k$ by $e^{\alpha_k}$, that is,
they are of the  form
\begin{align*}
& J_{r_k}(e^{\alpha_k}), \quad \alpha_k \in \reals, \\
& J_{r_k}(e^{\alpha_k})\quad\hbox{and}\quad  J_{r_k}(e^{\overline{\alpha}_k}),
\quad \alpha_k = \lambda_k + i\mu_k,\> \mu_k \not= 0.
\end{align*}
If  $\alpha_k \in \reals$, then $e^{\alpha_k} > 0$, so the
negative eigenvalues of $A$ must be of the form $e^{\alpha_k}$
or $e^{\overline{\alpha}_k}$, with  $\alpha_k$ complex.
This implies that
$\alpha_k = \lambda_k + (2h + 1)i\pi$, for some $h\in \integs$,
but then $\overline{\alpha}_k = \lambda_k - (2h + 1)i\pi$ and so
\[
e^{\alpha_k} = e^{\overline{\alpha}_k}.
\]
Consequently, negative eigenvalues of $A$ are associated with 
Jordan blocks that occur in pair, as claimed.
\end{proof}

\medskip
\remark
It can be shown (see Culver \cite{Culver})
that all the logarithms of a Jordan block,
$J_{r_k}(\alpha_k)$, corresponding to a real 
eigenvalue  $\alpha_k > 0$ are obtained by adding the matrices
\[
i 2\pi h_k I, \quad h_k \in \integs,
\]
to the solution given  by the proof of Theorem  \ref{log2} and
that all  the logarithms of a Jordan block,
$J_{2r_k}(\alpha_k, \beta_k)$, are obtained by adding the matrices
\[
i 2\pi h_k I + 2\pi l_k E
\quad h_k, l_k \in \integs,
\]
to the solution given  by the proof of Theorem  \ref{log2},
where
\[
E = 
\begin{pmatrix}
E_2 & \cdots & 0 \\
\vdots & \ddots & \vdots \\
0 & \cdots & E_2
\end{pmatrix},
\qquad
E_2 = 
\begin{pmatrix}
0 & -1 \\
1 & 0
\end{pmatrix}.
\]

\medskip
One should be careful no to relax the condition of Theorem \ref{log2}
to the more liberal
condition stating that for every Jordan block, $J_{r_k}(\alpha_k)$,
for which $\alpha_k < 0$, the dimension $r_k$ is even
({\it i.e\/}, $\alpha_k$ occurs an even number of times).
For example, the following matrix
\[
A = 
\begin{pmatrix}
-1 & 1 \\
0 & -1
\end{pmatrix}
\]
satisfies the more liberal condition but it does not possess
any real logarithm, as the reader will verify.
On the other hand, we have the following corollary:

\begin{cor}
\label{log3}
For every real invertible matrix, $A$, if $A$ has no negative
eigenvalues, then $A$ has a real logarithm.
\end{cor} 

\medskip
More results about the number of real logarithms of real matrices can be 
found in  Culver \cite{Culver}. In particular, Culver gives
a necessary and sufficient condition for a real matrix, $A$, to have
a unique real logarithm. This condition is quite strong. In particular,
it requires that all the eigenvalues of $A$ be real and positive.

\medskip
A different approach is to restrict the domain of real logarithms
to obtain a sufficient condition for the uniqueness of a logarithm.
We now discuss this approach. First, we state the following
property that will be useful later:

\begin{prop}
\label{multJod}
For every (real or complex) invertible matrix, $A$, there is 
a semisimple matrix, $S$, and a unipotent matrix, $U$, so that
\[
A = SU\qquad\hbox{and}\quad SU = US.
\]
Furthermore, $S$ and $U$ as above are unique.
\end{prop}

\begin{proof}
Proposition \ref{multJod} follows immediately from Theorem
\ref{Jordform1}, the details are left as an exercise.
\end{proof}

\medskip
The form, $SU$, of an invertible matrix is often called the 
{\it multiplicative Jordan decomposition\/}.

\begin{defin}
\label{stardom}
Let $\s{S}(n)$ denote the set of all  real matrices
whose  eigenvalues, $\lambda + i\mu$, lie in the horizontal strip
determined by the condition
$-\pi < \mu < \pi$.
\end{defin}

\medskip
It is easy to see that $\s{S}(n)$ is star-shaped (which means that if it 
contains $A$, then it contains $\lambda A$ for all $\lambda\in [0, 1]$)
and open (because the roots of a polynomial are continuous functions
of the coefficients of the polynomial).
As $\s{S}(n)$ is star-shaped, it is path-connected.
Furthermore, if $A\in \s{S}(n)$, then $PAP^{-1} \in \s{S}(n)$
for every invertible matrix, $P$. The remarkable property of
$\s{S}(n)$ is that the restriction of the exponential to $\s{S}(n)$
is a diffeomorphism onto its image.
To prove this fact we will need the following proposition:

\begin{prop}
\label{semiexp}
For any two real or complex
matrices, $S_1$ and $S_2$, if the eigenvalues, $\lambda + i\mu$, 
of $S_1$ and $S_2$ satisfy the condition 
$-\pi < \mu\leq \pi$, if 
$S_1$ and $S_2$ are semisimple and if $e^{S_1} = e^{S_2}$, then
$S_1 = S_2$.
\end{prop}

\begin{proof}
Since $S_1$ and $S_2$ are semisimple, they can be diagonalized
over $\complex$, so let
$(u_1, \ldots, u_n)$ be a basis of eigenvectors of $S_1$ 
associated with the (possibly complex) eigenvalues 
$\lambda_1, \ldots, \lambda_n$ and let
$(v_1, \ldots, v_n)$ be a basis of eigenvectors of $S_2$ 
associated with the (possibly complex) eigenvalues 
$\mu_1, \ldots, \mu_n$. We prove that if 
$e^{S_1} = e^{S_2} = A$, then
$S_1(v_i) = S_2(v_i)$
for all $v_i$, which shows that $S_1 = S_2$.

\medskip
Pick any eigenvector, $v_i$, of $S_2$ and write $v = v_i$ and $\mu = \mu_i$.
We have
\[
v = \alpha_1 u_1 + \cdots + \alpha_k u_k,
\]
for some unique $\alpha_j$'s. We compute $A(v)$ in two different ways.
We know that $e^{\mu_1}, \ldots, e^{\mu_n}$ are the eigenvalues of $e^{S_2}$
for the eigenvectors $v_1, \ldots, v_n$, so
\[
A(v) = e^{S_2}(v) = e^{\mu} v =
\alpha_1e^{\mu} u_1 + \cdots + \alpha_k e^{\mu} u_k.
\]
Similarly, we know that 
$e^{\lambda_1}, \ldots, e^{\lambda_n}$ are the eigenvalues of $e^{S_1}$
for the eigenvectors $u_1, \ldots, u_n$, so
\begin{eqnarray*} 
A(v) & = & A(\alpha_1 u_1 + \cdots + \alpha_k u_k) \\ 
 & = & \alpha_1 A(u_1) + \cdots + \alpha_k A(u_k) \\ 
 & = & \alpha_1 e^{S_1}(u_1) + \cdots + \alpha_k e^{S_1}(u_k) \\ 
 & = & \alpha_1 e^{\lambda_1} u_1 + \cdots + \alpha_k e^{\lambda_n} u_k.
\end{eqnarray*} 
Therefore, we deduce that
\[
\alpha_ke^{\mu} = \alpha_k e^{\lambda_k},
\qquad 1\leq k \leq n.
\]
Consequently, if $\alpha_k\not= 0$, then
\[
e^{\mu} = e^{\lambda_k},
\]
which implies $\mu - \lambda_k = i 2\pi h$,
for some $h\in \integs$.
However, due to the hypothesis on the eigenvalues of $S_1$ and  $S_2$,
$\mu$ and $\lambda_i$ must belong to the horizontal strip
determined by the condition
$-\pi < \Immag(z) \leq \pi$, so we must have $h = 0$ and then
$\mu = \lambda_k$.

\medskip
If we let 
$I = \{k \mid \lambda_k = \mu\}$ (which is nonempty since $v\not= 0$), 
then
$v = \sum_{k \in I} \alpha_k u_k$ and we have
\begin{eqnarray*} 
S_1(v) & = & S_1\left(\sum_{k \in I} \alpha_k u_k\right) \\
& = & \sum_{k \in I} \alpha_k S_1(u_k) \\
& = & \sum_{k \in I} \alpha_k \lambda_k u_k \\
& = & \sum_{k \in I} \alpha_k \mu u_k \\
& = & \mu \sum_{k \in I} \alpha_k  u_k  = \mu v.
\end{eqnarray*} 
Therefore, $S_1(v) = \mu v$. As $\mu$ is an eigenvector of $S_2$
for the eigenvalue $\mu$, we also have $S_2(v) = \mu v$.
Therefore,
\[
S_1(v_i) = S_2(v_i), \quad i = 1, \ldots, n,
\]
which proves that $S_1 = S_2$.
\end{proof}

\medskip
Obviously, Proposition \ref{semiexp} holds for real semisimple
matrices, $S_1, S_2$, in  $\s{S}(n)$, since the condition
for being in  $\s{S}(n)$ is $-\pi < \Immag(\alpha) < \pi$ 
for every eigenvalue, $\alpha$, of $S_1$ or $S_2$.

\medskip
We can now state our next theorem, an important result.
This theorem is a consequence of a more general fact proved
in Bourbaki \cite{BourbakiLie1}
(Chapter III, Section 6.9, Proposition 17, see also Theorem 6).

\begin{thm}
\label{stardiffeo}
The restriction of the exponential map to  $\s{S}(n)$ is 
a diffeomorphism of $\s{S}(n)$ onto its image, $\exp(\s{S}(n))$.
If $A\in \exp(\s{S}(n))$, then $PAP^{-1}\in \s{S}(n)$, for every (real)
invertible matrix, $P$.
Furthermore, $\exp(\s{S}(n))$ is an open subset of $\mathbf{GL}(n, \reals)$
containing $I$ and  $\exp(\s{S}(n))$ contains the open ball, 
$B(I, 1) = \{A \in \mathbf{GL}(n, \reals) \mid \norme{A - I} < 1\}$,
for every norm $\norme{\>}$ on $n\times n$ matrices
satisfying the condition $\norme{AB} \leq \norme{A}\norme{B}$.
\end{thm}

\begin{proof}
A complete proof is given in  Mmeimn\'e and Testard \cite{Mneimne},
Chapter 3, Theorem 3.8.4.
Part of the proof consists in showing that $\exp$ is a local 
diffeomorphism and for this, to prove that
$d\exp(X)$ is invertible. This requires finding
an explicit formula for the derivative of the exponential
and we prefer to omit this computation, which is quite
technical. Proving that $B(I, 1) \subseteq \s{S}(n)$ is easier
but requires a little bit of complex analysis.
Once these facts are established, it remains to prove that
$\exp$ is injective on $\s{S}(n)$, which we will prove.

\medskip
The trick is to use both the Jordan decomposition and the
multiplicative Jordan decomposition!
Assume that $X_1, X_2\in \s{S}(n)$ and that $e^{X_1} = e^{X_2}$. 
Using Theorem \ref{Jordform1}
we can write $X_1 = S_1 + N_1$ and $X_2 = S_2 + N_2$,
where $S_1, S_2$ are semisimple, $N_1, N_2$ are nilpotent,
$S_1N_1 = N_1S_1$,   and $S_2N_2 = N_2S_2$.
From  $e^{X_1} = e^{X_2}$, we get
\[
e^{S_1}e^{N_1} = e^{S_1 + N_1} =  e^{S_2 + N_2} = e^{S_2}e^{N_2}.
\]
Now, $S_1$  and $S_2$ are semisimple, so $e^{S_1}$ and $e^{S_2}$
are semisimple and $N_1$ and $N_2$ are nilpotent so
$e^{N_1}$ and $e^{N_2}$ are unipotent. Moreover,
as $S_1N_1 = N_1S_1$ and $S_2N_2 = N_2S_2$, 
we have $e^{S_1}e^{N_1} = e^{N_1}e^{S_1}$ and
$e^{S_2}e^{N_2} = e^{N_2}e^{S_2}$.
By the uniqueness property of Proposition \ref{multJod},
we conclude that
\[
e^{S_1} = e^{S_2}
\quad\hbox{and}\quad
e^{N_1} = e^{N_2}.
\]
Now, as $N_1$ and $N_2$ are nilpotent, there is some $r$ so that 
$N_1^r = N_2^r = 0$ and then, it is clear that
$e^{N_1} = I + \widetilde{N}_1$ and $e^{N_2} = I + \widetilde{N}_2$
with $\widetilde{N}_1^r = 0$ and   $\widetilde{N}_2^r = 0$.
Therefore, we can apply Proposition \ref{expnil} to conclude that
\[
N_1 = N_2.
\]

As $S_1, S_2\in \s{S}(n)$ are semisimple and
$e^{S_1}  = e^{S_2}$, by Proposition \ref{semiexp}, we conclude that
\[
S_1 = S_2.
\]
Therefore, we finally proved that $X_1 = X_2$,
showing that $\exp$ is injective on $\s{S}(n)$.
\end{proof}

\medskip
\remark
Since proposition \ref{semiexp} holds for semisimple
matrices, $S$, such that the condition
$-\pi < \mu\leq \pi$ holds for every eigenvalue, $\lambda + i\mu$, 
of $S$, the restriction of the exponential to real matrices, $X$,
whose eigenvalues satisfy this condition is  injective.
Note that the image of these matrices under the exponential contains
matrices, $A = e^X$, with negative eigenvalues.
Thus, combining Theorem \ref{log2} and the above injectivity
result we could state an existence and uniqueness result for 
real logarithms of real matrices that is more general than
Theorem \ref{logexistuniq} below.
However this is not a practical result since
it requires a condition on the number of Jordan blocks 
and such a condition is hard to check. Thus,
we will restrict ourselves to real matrices with no negative
eigenvalues (see Theorem \ref{logexistuniq}).

\medskip
Since the eigenvalues of a nilpotent matrix are zero and since
symmetric matrices have real eigenvalues, Theorem \ref{stardiffeo}
has has two interesting corollaries.
Denote by $\mathbf{S}(n)$ the vector space of real $n\times n$ matrices
and by $\mathbf{SPD}(n)$ the set of $n\times n$ symmetric, positive, definite
matrices. It is known that 
$\mapdef{\exp}{\mathbf{S}(n)}{\mathbf{SPD}(n)}$ is a bijection.

\begin{cor}
\label{stardiffeocor}
The exponential map has the following properties:
\begin{enumerate}
\item[(1)]
The map
$\mapdef{\exp}{\s{N}\mathit{il}(r)}{\s{U}\mathit{ni}(r)}$ is a diffeomorphism.
\item[(2)]
The map
$\mapdef{\exp}{\mathbf{S}(n)}{\mathbf{SPD}(n)}$ is a diffeomorphism.
\end{enumerate}
\end{cor}

\medskip
By combining Theorem \ref{log2} and Theorem \ref{stardiffeo} we obtain
the following result about the existence and uniqueness of logarithms of
real matrices:

\begin{thm}
\label{logexistuniq}
(a)
If $A$ is any real invertible $n\times n$
matrix and $A$ has no negative eigenvalues,
then  $A$ has a unique real logarithm, $X$, with $X\in \s{S}(n)$.

\medskip
(b)
The image, $\exp(\s{S}(n))$, of $\s{S}(n)$ by the exponential map
is the set of real invertible matrices with no negative eigenvalues
and $\mapdef{\exp}{\s{S}(n)}{\exp(\s{S}(n))}$ is a diffeomorphism
between these two spaces.
\end{thm}

\begin{proof}
(a)
If we go back to the proof of Theorem  \ref{log2}, we see that complex
eigenvalues of the logarithm, $X$, produced by that proof only occur
for matrices
\[
S(\rho_k, \theta_k) = 
\begin{pmatrix}
\log\rho_k & - \theta_k \\
\theta_k & \log\rho_k
\end{pmatrix},
\]
associated with eigenvalues 
$\lambda_k + i \mu_k = \rho_k\,e^{i\theta_k}$.
However, the eigenvalues of such matrices are
$\log\rho_k \pm i\theta_k$ and since $A$ has no negative eigenvalues, we may
assume that $-\pi < \theta_k < \pi$, and so 
$X \in \s{S}(n)$, as desired. By Theorem  \ref{stardiffeo},
such a logarithm is unique.

\medskip
(b)
Part (a) proves that the set of real invertible matrices with no 
negative eigenvalues is contained in $\exp(\s{S}(n))$.
However, for any matrix, $X\in \s{S}(n)$, 
since every eigenvalue of $e^X$ is of the form
$e^{\lambda + i\mu} = e^{\lambda}e^{i\mu}$ for some eigenvalue,
$\lambda + i\mu$, of $X$  and since $\lambda + i\mu$ satisfies the condition 
$-\pi < \mu < \pi$, the number, $e^{i\mu}$, is never negative, 
so $e^X$  has no negative eigenvalues.
Then, (b) follows directly from Theorem \ref{stardiffeo}.
\end{proof}

\medskip
\remark
Theorem \ref{logexistuniq} (a) first appeared in
Kenney and Laub \cite{KenneyLaub} (Lemma A2, Appendix A) but
without proof.

\section{Square Roots of Real Matrices; 
Criteria for Existence and Uniqueness}
\label{sec3}
In this section we investigate the problem of finding a square root
of a matrix, $A$, that is, a matrix, $X$, such that $X^2 = A$.
If $A$ is an invertible (complex) matrix, then it always has
a square  root, but singular matrices may fail to have a square root.
For example, the nilpotent matrix, 
\[
H = 
\begin{pmatrix}
0 & 1 & 0 &  \cdots &  0 \\
0 & 0  & 1  & \cdots & 0 \\
\vdots & \vdots  & \ddots  & \ddots  &  \vdots      \\
0 & 0 & 0 & \ddots &  1 \\
0 & 0 & 0 & \cdots &  0
\end{pmatrix}
\] 
has no square root (ckeck this!). The problem of finding square roots
of matrices is thoroughly investigated in Gantmacher \cite{Gantmacher1}, 
Chapter VIII, Sections 6 and 7. For singular matrices, finding a square root
reduces to the problem of finding the square root of a nilpotent matrix,
which is not always possible.
A necessary and sufficient condition
for the existence of a square root
is given in Horn and Johnson \cite{HornJohn2}, see Chapter 6, Section 4,
especially Theorem 6.1.12 and Theorem 6.4.14.
This criterion is rather complicated because  
its deals with non-singular as well
as singular matrices. In this paper,
we will restrict our attention to invertible matrices.
The main two Theorems of this section are Theorem \ref{sqmat2} and
Theorem \ref{squareexistuniq}. The former theorem appears in
Higham \cite{Higham87} (Theorem 5).
The first step is to prove a version of Theorem \ref{expJord}
for the function $A \mapsto A^2$, where $A$ is invertible.
In this section, we state the following theorem without proof.
A proof is given in Section \ref{sec6}.

\begin{thm}
\label{squareJord}
For any (real or complex) invertible $n\times n$ matrix, $A$, if
$A = P J P^{-1}$ where $J$ is a Jordan matrix of the form
\[
J = 
\begin{pmatrix}
J_{r_1}(\lambda_1) &  \cdots & 0 \\
 \vdots  & \ddots   & \vdots    \\
 0 & \cdots & J_{r_m}(\lambda_m)
\end{pmatrix},
\]
then there is some invertible matrix, $Q$, so that
the  Jordan form of $A^2$ is given by
\[
e^A = Q\, s(J)\,  Q^{-1},
\]
where $s(J)$ is the Jordan matrix
\[
s(J) = 
\begin{pmatrix}
J_{r_1}(\lambda_1^2) &  \cdots & 0 \\
 \vdots  & \ddots   & \vdots    \\
 0 & \cdots & J_{r_m}(\lambda_m^2)
\end{pmatrix},
\]
that is, each $J_{r_k}(\lambda_k^2)$ is obtained from $J_{r_k}(\lambda_k)$
by replacing all the diagonal enties $\lambda_k$ by $\lambda_k^2$.
Equivalently, if the list of elementary divisors of $A$ is
\[
(X - \lambda_1)^{r_1}, \ldots, (X - \lambda_m)^{r_m},
\]
then the list of elementary divisors of $A^2$ is
\[
(X - \lambda_1^2)^{r_1}, \ldots, (X - \lambda_m^2)^{r_m}.
\]
\end{thm}

\medskip

\remark
Theorem \ref{squareJord} can be easily generalized to the map
$A \mapsto A^p$, for any $p \geq 2$, that is, by replacing $A^2$ by $A^p$,
provided $A$ is invertible. Thus, if the list of elementary divisors of $A$ is
\[
(X - \lambda_1)^{r_1}, \ldots, (X - \lambda_m)^{r_m},
\]
then the list of elementary divisors of $A^p$ is
\[
(X - \lambda_1^p)^{r_1}, \ldots, (X - \lambda_m^p)^{r_m}.
\]

\medskip
The next step is to find the square root of a Jordan block.
Since we are assuming that our matrix is invertible, every Jordan block,
$J_{r_k}(\alpha_k)$, can be written as
\[
J_{r_k}(\alpha_k) = \alpha_k I\left(I + \frac{H}{\alpha_k} \right),
\]
where $H$ is nilpotent.
It is easy to find a square root of $\alpha_k I$. If 
$\alpha_k = \rho_k e^{\theta_k}$, with $\rho_k > 0$, then
\[
S_k = 
\begin{pmatrix}
\sqrt{\rho_k}\, e^{i\frac{\theta_k}{2}} & 0 & \cdots & 0 \\
0 & \sqrt{\rho_k}\, e^{i\frac{\theta_k}{2}}  &  \cdots & 0 \\
\vdots & \vdots  & \ddots   &  \vdots      \\
0 & 0 &  \cdots &  \sqrt{\rho_k}\, e^{i\frac{\theta_k}{2}}
\end{pmatrix}
\]
is a square root of $\alpha_k I$.
Therefore,  the problem reduces
to finding  square roots of unipotent matrices. For this, we recall the
power series
\begin{eqnarray*}
(1 + x)^{\frac{1}{2}} & = & 1 + \frac{1}{2} x + \cdots + 
\frac{1}{n!}\frac{1}{2}\left(\frac{1}{2} - 1\right)  \cdots 
\left(\frac{1}{2} - n + 1\right) x^{n} + \cdots \\
& = & \sum_{n = 0}^{\infty} (-1)^{n - 1} \frac{(2n)!}{(2n - 1)(n!)^2 2^{2n}}\, x^n, 
\end{eqnarray*}
which is normally convergent for $|x | < 1$.
Then, we can define the power series, $R$, of a matrix variable, $A$, by
\[
R(A) = 
\sum_{n = 1}^{\infty} (-1)^{n - 1} \frac{(2n)!}{(2n - 1)(n!)^2 2^{2n}}\,
 A^n, 
\]
and this power series  converges normally for $\norme{A} < 1$.
As a formal power series, note that $R(0) = 0$ 
and $R'(0) = \frac{1}{2}  \not= 0$
so, by a theorem about formal power series, $R$ has a unique inverse, $S$,
such that $S(0) = 0$ (see Lang \cite{Langcomplex} or H. Cartan
\cite{Cartancomplex}).
But, if we consider the power series, $S(A) = (I + A)^2 - I$, when $A$ is
a real number, we have $R(A) = \sqrt{1 + A} - 1$, so we get
\[
R\circ S(A) = \sqrt{1 + (1 + A)^2 - 1} - 1 = A,
\]
from wich we deduce that
$S$ and $R$ are mutual inverses. But, $R$ converges everywhere and $S$ converges
for $\norme{A} < 1$, so by another theorem about converging power series,
if we let $\sqrt{I + A} = R(A) + I$, there is some $r$, 
with $0 < r < 1$, so that
\[
(\sqrt{I + A})^2 = I + A, \qquad\hbox{if}\quad  \norme{A} < r 
\]
and 
\[
\sqrt{(I + A)^2} = I + A,  \qquad\hbox{if}\quad  \norme{A} <  r.
\]
If $A$ is unipotent, that is, $A = I + N$ with $N$ nilpotent, we see that
the series has only finitely many terms. This fact allows us to prove
the proposition below.

\begin{prop}
\label{squareunip}
The squaring map, $A \mapsto A^2$, is a homeomorphism from
$\s{U}\mathit{ni}(r)$ to itself whose inverse is the map 
$A \mapsto \sqrt{A} = R(A - I) + I$.
\end{prop}

\begin{proof}
If $A = I + N$ with $N^r = 0$, 
as $A^2 = I + 2N + N^2$ it is clear that $(2N + N^2)^r = 0$,
so the squaring map is well defined on unipotent matrices. 
We use the technique of Proposition \ref{expnil}.
Consider the map 
\[
t \mapsto (\sqrt{I + t N})^2 - (I + t N), \qquad t\in \reals.
\]
It is a polynomial since $N^r = 0$. Furthermore, for $t$ sufficiently small,
$\norme{t N} < 1$ and we have $(\sqrt{I + t N})^2 = I + t N$,
so the above polynomial
vanishes in a neighborhood of $0$, which implies that it is 
identically zero. Therefore,  $(\sqrt{I +  N})^2 = I +  N$,
as required.

\medskip
Next, consider the map 
\[
t \mapsto \sqrt{(I + t N)^2} - (I + t N), \qquad t\in \reals.
\]
It is a polynomial since $N^r = 0$. Furthermore, for $t$ sufficiently small,
$\norme{t N} < 1$ and we have $\sqrt{(I + t N)^2} = I + t N$,
so we conclude as above that the above map is identically zero and that
$\sqrt{(I +  N)^2} = I +  N$.
\end{proof}

\medskip
\remark
Proposition \ref{squareunip} can be easily generalized to the map
$A \mapsto A^p$, for any $p \geq 2$, by using the power series
\[
(I + A)^{\frac{1}{p}} = 
 I + \frac{1}{p} A + \cdots + 
\frac{1}{n!}\frac{1}{p}\left(\frac{1}{p} - 1\right)  \cdots 
\left(\frac{1}{p} - n + 1\right) A^{n} + \cdots .
\]

\medskip
Using proposition \ref{squareunip}, we can find a square root
for the unipotent part of a Jordan block,
\[
J_{r_k}(\alpha_k) = \alpha_k I\left(I + \frac{H}{\alpha_k} \right).
\]
If  $N_k = \frac{H}{\alpha_k}$, then 
\[
\sqrt{I + N_k} = I +  
\sum_{j = 1}^{r_k - 1} (-1)^{j - 1} \frac{(2j)!}{(2j - 1)(j!)^2 2^{2j}}\,  N_k^j 
\]
is a square root of $I + N_k$.
Therefore, we obtained the following theorem:

\begin{thm}
\label{sqmat1}
Every (complex) invertible matrix, $A$, has a square root.
\end{thm}

\medskip
\remark
Theorem \ref{sqmat1} can be easily generalized to $p^{\mathrm{th}}$ roots,
for any $p \geq 2$,

\medskip
We now consider the problem of finding a real square root of 
an invertible real matrix.
It turns out that the necessary and sufficient condition is exactly
the condition for finding a real logarithm of a real matrix.

\begin{thm}
\label{sqmat2}
Let $A$ be a real invertible $n\times n$ matrix and let 
$(X - \alpha_1)^{r_1}, \ldots, (X - \alpha_m)^{m_1}$
be its list of elementary divisors or, equivalently,
let $J_{r_1}(\alpha_1), \ldots$, $J_{r_m}(\alpha_m)$
be its list of Jordan blocks. Then, $A$ has a real square root
iff for every $r_i$ and every real eigenvalue $\alpha_i < 0$, 
the number, $m_i$,  of Jordan blocks identical to $J_{r_i}(\alpha_i)$
is even.
\end{thm}

\begin{proof}
The proof is very similar to the proof of Theorem \ref{log2} so 
we only point out the necessary changes.
Let $J'$ be  a real Jordan matrix so that 
\[
A = P J' P^{-1},
\]
where $J'$ satisfies  conditions (1) and (2) of Theorem  \ref{Jordanform3}.
As $A$ is invertible, every block of $J'$ of the form $J_{r_k}(\alpha_k)$
corresponds to a real eigenvalue with  $\alpha_k > 0$ and we can write
$J_{r_k}(\alpha_j) = \alpha_k I(I + N_k)$, where $N_k$ is nilpotent.
As in Theorem \ref{sqmat1}, we can find a real square root, $M_k$, of $I + N_k$
and as  $\alpha_k > 0$, the diagonal matrix  $\alpha_k I$
has the real square root
\[
S_k = 
\begin{pmatrix}
\sqrt{\alpha_k} & 0 & \cdots & 0 \\
0 &  \sqrt{\alpha_k}   &  \cdots & 0 \\
\vdots & \vdots  & \ddots   &  \vdots      \\
0 & 0 &  \cdots &  \sqrt{\alpha_k} 
\end{pmatrix}.
\]
Set $Y_k = S_kM_k$.

\medskip
The other real Jordan blocks of $J'$ are of the form
$J_{2r_k}(\lambda_k, \mu_k)$,
with $\lambda_k, \mu_k \in \reals$, not both zero. Consequently,
we can write
\[
J_{2r_k}(\lambda_k, \mu_k) =   D_k(I + N_k) 
\]
where 
\[
D_k = 
\begin{pmatrix}
L(\lambda_k, \mu_k) &  \cdots & 0 \\
 \vdots  & \ddots   & \vdots    \\
 0 & \cdots & L(\lambda_k, \mu_k)
\end{pmatrix}
\]
with
\[
L(\lambda_k, \mu_k) = 
\begin{pmatrix}
\lambda_k & - \mu_k \\
\mu_k & \lambda_k
\end{pmatrix},
\]
and  $N_k = D_k^{-1}H_k$ is nilpotent. 
We can find a square root, $M_k$, of $I + N_k$ as in Theorem \ref{sqmat1}.
If we write $\lambda_k + i\mu_k =  \rho_k e^{i \theta_k}$, then
\[
L(\lambda_k, \mu_k) = 
\rho_k
\begin{pmatrix}
\cos\theta_k & - \sin\theta_k \\
\sin\theta_k & \cos\theta_k
\end{pmatrix}.
\] 
Then, if we set
\[
S(\rho_k, \theta_k) = 
\sqrt{\rho_k}
\begin{pmatrix}
\cos\left(\frac{\theta_k}{2}\right) & - \sin\left(\frac{\theta_k}{2}\right) \\
\sin\left(\frac{\theta_k}{2}\right) & \cos\left(\frac{\theta_k}{2}\right)
\end{pmatrix},
\]
a real matrix, we have
\[
L(\lambda_k, \mu_k) = S(\rho_k, \theta_k)^2.
\]
If we form the real block diagonal matrix,
\[
S_k = 
\begin{pmatrix}
S(\rho_k, \theta_k) &  \cdots & 0 \\
 \vdots  & \ddots   & \vdots    \\
 0 & \cdots & S(\rho_k, \theta_k)
\end{pmatrix},
\]
we have $D_k = S_k^2$ and then
the matrix $Y_k = S_kM_k$ is a square root
of $J_{2r_k}(\lambda_k, \mu_k)$.
Finally, if $Y$ is the block diagonal matrix
$\mathrm{diag}(Y_1, \ldots, Y_m)$, then
$X = P Y P^{-1}$ is a square root of $A$.

\medskip
Let us now prove that if $A$ has a real square root, $X$, then $A$ satisfies
the condition of Theorem \ref{sqmat2}. 
Since $X$ is a real matrix,
we know from the proof of Theorem \ref{Jordanform2} that
the Jordan blocks of $X$ associated with complex eigenvalues occur
in conjugate pairs, so they are of the form
\begin{align*}
& J_{r_k}(\alpha_k), \quad \alpha_k \in \reals, \\
& J_{r_k}(\alpha_k)\quad\hbox{and}\quad  J_{r_k}(\overline{\alpha}_k),
\quad \alpha_k = \lambda_k + i\mu_k,\> \mu_k \not= 0.
\end{align*}
By Theorem \ref{squareJord}, the Jordan blocks of $A = X^2$ are 
obtained by replacing each $\alpha_k$ by $\alpha_k^2$, that is,
they are of the  form
\begin{align*}
& J_{r_k}(\alpha_k^2), \quad \alpha_k \in \reals, \\
& J_{r_k}(\alpha_k^2)\quad\hbox{and}\quad  J_{r_k}(\overline{\alpha}_k^2),
\quad \alpha_k = \lambda_k + i\mu_k,\> \mu_k \not= 0.
\end{align*}
If  $\alpha_k \in \reals$, then $\alpha_k^2 > 0$, so the
negative eigenvalues of $A$ must be of the form $\alpha_k^2$
or $\overline{\alpha}_k^2$, with  $\alpha_k$ complex.
This implies that
$\alpha_k = \sqrt{\rho_k}\, e^{i\frac{\pi}{2}}$, 
but then $\overline{\alpha}_k = \sqrt{\rho_k}\, e^{-i\frac{\pi}{2}}$ and so
\[
\alpha_k^2 = \overline{\alpha}_k^2.
\]
Consequently, negative eigenvalues of $A$ are associated with 
Jordan blocks that occur in pair, as claimed.
\end{proof}

\medskip
\remark
Theorem \ref{sqmat2} can be easily generalized to $p^{\mathrm{th}}$ roots,
for any $p \geq 2$,

\medskip
Theorem \ref{sqmat2} appears in  Higham \cite{Higham87} as Theorem 5
but no explicit proof is given. Instead, Higham states:
``The proof is a straightfoward modification of Theorem 1 in
Culver  \cite{Culver} and is omitted.'' Culver's proof
uses results from  Gantmacher \cite{Gantmacher1} and 
does not provide a constructive method for obtaining a square root. 
We gave a more constructive proof (but perhaps longer).

\begin{cor}
\label{sqmat3}
For every real invertible matrix, $A$, if $A$ has no negative
eingenvalues, then $A$ has a real square root.
\end{cor}

\medskip
We will now provide a sufficient condition for the uniqueness
of a real square root. For this, we consider the open set, $\s{H}(n)$,
consisting of all real $n\times n$  matrices whose eigenvalues,
$\alpha = \lambda + i\mu$, have a positive real part, $\lambda > 0$.
We express this condition as $\Re(\alpha) > 0$.
Obviously, such matrices are invertible  and can't have negative eigenvalues.
We need a version of Proposition \ref{semiexp} for semisimple matrices
in $\s{H}(n)$.

\medskip
\remark
To deal with $p^{\mathrm{th}}$ roots, we consider matrices whose
eigenvalues, $\rho e^{i \theta}$, satisfy the condition
$-\frac{\pi}{p} < \theta < \frac{\pi}{p}$.

\begin{prop}
\label{semisquare}
For any two real or complex matrices, $S_1$ and  $S_2$, 
if the eigenvalues, $\rho e^{i\theta}$, of $S_1$ and $S_2$
satisfy the condition $-\frac{\pi}{2} < \theta \leq  \frac{\pi}{2}$,
if $S_1$ and $S_2$ are semisimple and if $S_1^2 = S_2^2$, then
$S_1 = S_2$.
\end{prop}

\begin{proof}
The proof is very similar to that of Proposition \ref{semiexp}
so we only indicate where  modifications are needed.
We use the fact that if $u$ is an eigenvector of a linear map, $A$,
associated with some eigenvalue, $\lambda$, then 
$u$ is an eigenvector of  $A^2$
associated with the eigenvalue $\lambda^2$. We replace
every occurrence of $e^{\lambda_i}$ by $\lambda_i^2$
(and $e^{\mu}$ by $\mu^2$).
As in the proof of Proposition \ref{semiexp}, we obtain the equation
\[
\alpha_1\mu^2 u_1 + \cdots + \alpha_k \mu^2 u_k
= \alpha_1 \lambda_1^2 u_1 + \cdots + \alpha_k \lambda^2 u_k.
\]
Therefore, we deduce that
\[
\alpha_k \mu^2 = \alpha_k \lambda_k^2,
\qquad 1\leq k \leq n.
\]
Consequently, as $\mu, \lambda_k \not= 0$, if $\alpha_k\not= 0$, then
\[
\mu^2 = \lambda_k^2,
\]
which implies $\mu =  \pm\lambda_k$.
However, the hypothesis on the eigenvalues of  
$S_1$ and  $S_2$ implies that $\mu = \lambda_k$.
The end of the proof is identical to that of Proposition
\ref{semiexp}.
\end{proof}

\medskip
Obviously, Proposition \ref{semisquare} holds for real semisimple
matrices, $S_1, S_2$, in  $\s{H}(n)$.

\medskip
\remark
Proposition \ref{semisquare} also holds
for the map $S \mapsto S^p$,
for any $p \geq 2$, under the condition
$-\frac{\pi}{p} < \theta \leq  \frac{\pi}{p}$.

\medskip
We have the following analog of Theorem \ref{stardiffeo},
but we content ourselves with a weaker result:

\begin{thm}
\label{stardiffeo2}
The restriction of the squaring map, $A \mapsto A^2$,  to $\s{H}(n)$ is 
injective.
\end{thm}

\begin{proof}
Let $X_1, X_2\in \s{H}(n)$ and assume that $X_1^2 = X_2^2$.
As $X_1$ and $X_2$ are invertible, by Proposition \ref{multJod},
we can write
$X_1 = S_1(I + N_1)$  and $X_2 = S_2(I + N_2)$, where
$S_1, S_2$ are semisimple, $N_1, N_2$ are nilpotent,
$S_1(I + N_1) = (I + N_1)S_1$ and $S_2(I + N_2) = (I + N_2)S_2$.
As $X_1^2 = X_2^2$, we get
\[
S_1^2(I + N_1)^2 = S_2^2(I + N_2)^2.
\]
Now, as $S_1$ and $S_2$ are semisimple and
invertible,  $S_1^2$ and $S_2^2$ are semisimple and invertible, 
and as $N_1$ and $N_2$ are nilpotent, $2N_1 + N_1^2$ and 
$2N_2 + N_2^2$ are nilpotent, so $(I + N_1)^2$ and $(I + N_2)^2$
are unipotent. Moreover, 
$S_1(I + N_1) = (I + N_1)S_1$ and $S_2(I + N_2) = (I + N_2)S_2$ imply that
$S_1^2(I + N_1)^2 = (I + N_1)^2S_1^2$ and $S_2^2(I + N_2)^2 = (I + N_2)^2S_2^2$.
Therefore, by the uniqueness statement of Proposition \ref{multJod},
we get
\[
S_1^2 = S_2^2
\quad\hbox{and}\quad
 (I + N_1)^2 = (I + N_2)^2.
\]
However, as  $X_1, X_2\in \s{H}(n)$ we have
$S_1, S_2\in \s{H}(n)$ and Proposition \ref{semisquare}
implies that $S_1 = S_2$. Since $I + N_1$ and $I + N_2$ are unipotent,
proposition \ref{squareunip} implies that $N_1 = N_2$.
Therefore, $X_1 = X_2$, as required.
\end{proof}

\medskip
\remark
Theorem \ref{stardiffeo2} also holds for the restriction of the
squaring map to real or complex matrices, $X$,
whose eigenvalues, $\rho e^{i\theta}$, satisfy the condition
$-\frac{\pi}{2} < \theta \leq  \frac{\pi}{2}$.
This result is proved in DePrima and Johnson \cite{DePrimaJohnson} 
by a different method. However, DePrima and Johnson
need an extra condition, see the discussion at the end of this section.

\medskip
We can now prove the analog of Theorem \ref{logexistuniq}
for square roots.

\begin{thm}
\label{squareexistuniq}
If $A$ is any real invertible $n\times n$ matrix and $A$ 
has no negative eigenvalues,
then  $A$ has a unique real  square root, $X$, with $X\in \s{H}(n)$.
\end{thm}

\begin{proof}
If we go back to the proof of Theorem  \ref{sqmat2}, we see that complex
eigenvalues of the square root, $X$, produced by that proof only occur
for matrices
\[
S(\rho_k, \theta_k) = 
\sqrt{\rho_k}
\begin{pmatrix}
\cos\left(\frac{\theta_k}{2}\right) & - \sin\left(\frac{\theta_k}{2}\right) \\
\sin\left(\frac{\theta_k}{2}\right) & \cos\left(\frac{\theta_k}{2}\right)
\end{pmatrix},
\]
associated with eigenvalues 
$\lambda_k + i \mu_k = \rho_k\, e^{i \theta_k}$.
However, the eigenvalues of such matrices are
$\sqrt{\rho_k}\, e^{\pm i \frac{\theta_k}{2}}$ 
and since $A$ has no negative eigenvalues, we may
assume that $-\pi < \theta_k < \pi$, and so 
$-\frac{\pi}{2} < \frac{\theta_k}{2} < \frac{\pi}{2}$, wich means that
$X \in \s{H}(n)$, as desired. By Theorem  \ref{stardiffeo2},
such a square root is unique.
\end{proof}

\medskip
Theorem \ref{squareexistuniq} is stated  in a number of
papers including Bini, Higham and Meini \cite{BiniHigham},
Cheng,  Higham, Kenney and  Laub \cite{Higham01} and 
Kenney and Laub \cite{KenneyLaub}.
Theorem \ref{squareexistuniq} also appears in Higham \cite{HighamFM}
as Theorem 1.29. Its proof relies on Theorem 1.26 and
Theorem 1.18 (both in Higham's book), whose proof is not given in full
(closer examination reveals that Theorem 1.36 (in Higham's book)
is needed to prove Theorem 1.26).
Although Higham's Theorem 1.26 implies our Theorem \ref{stardiffeo2}
we feel that the proof of Theorem \ref{stardiffeo2} is of
independent interest and is more direct.

\medskip
As we already said in Section \ref{sec2}, Kenney and Laub  \cite{KenneyLaub}
state Theorem  \ref{squareexistuniq} as Lemma A1 in Appendix A.
The proof is sketched  briefly. Existence follows from
the Cauchy integral formula for operators, a method used by
DePrima and Johnson \cite{DePrimaJohnson} in which a similar result is proved
for complex matrices (Section 4, Lemma 1). Uniqueness is proved in
DePrima and Johnson \cite{DePrimaJohnson} but it uses an extra condition.  
The hypotheses of Lemma 1 in 
DePrima and Johnson are that $A$ and $X$ are {\it complex\/} invertible 
matrices and that $X$  satisfies the conditions
\begin{enumerate}
\item[(i)]
$X^2 = A$,
\item[(ii)]
the eigenvalues, $\rho\, e^{i\theta}$, of $X$ satisfy
$-\frac{\pi}{2} < \theta \leq \frac{\pi}{2}$,
\item[(iii)]
For any matrix, $S$, if $AS = SA$, then $XS = SX$.
\end{enumerate}

\medskip
Observe that condition (ii) allows $\theta =  \frac{\pi}{2}$,
which yields matrices, $A = X^2$, with negative eigenvalues.
In this case, $A$ may not have any real square root but
DePrima and Johnson are only concerned with {\it complex\/} matrices
and a complex square root always exists.
To guarantee the existence of real logarithms, Kenney and Laub
tighten condition (ii) to  $-\frac{\pi}{2} < \theta < \frac{\pi}{2}$.
They also assert that
condition (iii) follows from conditions (i) and (ii).
This can be shown as follows: First, recall that
we have shown that uniqueness follows from (i) and (ii). 
Uniqueness under conditions (i) and (ii)
can also be shown to be a consequence of
Theorem 2 in Higham \cite{Higham87}.
Now, assume $X^2 = A$ and $SA = SA$. We may assume that $S$ is invertible
since the set of invertible matrices is dense in the set of all matrices.
Then, as $SA = AS$, we have
\[
(SXS^{-1})^2 = SX^2S^{-1} = SAS^{-1} = A.
\]
Thus, $SXS^{-1}$ is a square root of $A$. Furthermore, $X$ and $SXS^{-1}$
have the same eigenvalues so  $SXS^{-1}$ satisfies (i) and (ii) and, 
by uniqueness,
$X = SXS^{-1}$, that is, $XS = SX$.

\medskip
Since Kenney and Laub only provide a sketch of Theorem A1
and since Higham \cite{HighamFM} does not give all the details of the proof
either,
we felt that the reader would appreciate seeing
a complete proof of Theorem  \ref{squareexistuniq}.

\section{Conclusion}
\label{sec5}
It is interesting that Theorem \ref{logexistuniq} and Theorem
\ref{squareexistuniq} are the basis for numerical methods
for computing the exponential or the logarithm of a matrix.
The key point is that the following identities hold:
\[
e^{A} = (e^{A/2^k})^{2^k} 
\quad\hbox{and}\quad
\log(A) = 2^k \log(A^{1/2^k}),
\]
where in the second case, $A^{1/2^k}$ is the unique
$k$th square root of $A$ whose eigenvalues, $\rho\, e^{i\theta}$, lie in the sector
$-\frac{\pi}{2^k} < \theta  < \frac{\pi}{2^k}$.
The first identity is trivial and the second one can be shown by induction
from the identity
\[
\log(A) = 2 \log(A^{1/2}),
\]
where $A^{1/2}$ is the unique
square root of $A$ whose eigenvalues, $\rho\, e^{i\theta}$, lie in the sector \\
$-\frac{\pi}{2} < \theta  < \frac{\pi}{2}$.
Let $\widetilde{X} = A^{1/2}$, whose eigenvalues, $\rho\, e^{i\theta}$, 
lie in the sector 
$-\frac{\pi}{2} < \theta  < \frac{\pi}{2}$. Then, it is easy to see
that the eigenvalues, $\alpha$,
of $\log(\widetilde{X})$ satisfy the condition
$-\frac{\pi}{2} < \Immag(\alpha) < \frac{\pi}{2}$. Then,
$X = 2\log(\widetilde{X}) = 2 \log(A^{1/2})$ satisfies 
\[
e^X = e^{\log(A^{1/2}) + \log(A^{1/2})} = e^{\log(A^{1/2})}e^{\log(A^{1/2})} =
A^{1/2} A^{1/2} = A, 
\]
and the eigenvalues, $\alpha$,
of $X$ satisfy the condition $-\pi < \Immag(\alpha) < \pi$ so,
by the uniqueness part of Theorem \ref{logexistuniq}, we must have
$\log(A) =  2 \log(A^{1/2})$.

\medskip
The identity $\log(A) = 2^k \log(A^{1/2^k})$ leads to a
numerical method for computing the logarithm of a (real) matrix
first introduced by Kenney and Laub
known as the {\it inverse scaling and squaring algorithm\/},
see Kenney and Laub \cite{KenneyLaub}
and Cheng,  Higham, Kenney and  Laub \cite{Higham01}.
The idea is that if $A$ is close to the identity, then
$\log(A)$ can be computed accurately using either a truncated
power series expansion of $\log(A)$ or better, rational approximations
know as {\it Pad\'e approximants\/}. In order to bring $A$ close to the 
identity, iterate the operation of taking the square root of $A$ to obtain
$A^{1/2^k}$. Then, after having computed $\log(A^{1/2^k})$,
scale $\log(A^{1/2^k})$ by the factor $2^k$. For details of this method,
see  Kenney and Laub \cite{KenneyLaub}
and Cheng,  Higham, Kenney and  Laub \cite{Higham01}.
The inverse squaring and scaling method plays an important role in the
{\it log-Euclidean framework\/} introduced by
Arsigny, Fillard, Pennec and Ayache, see Arsigny \cite{Arsignythesis}, 
Arsigny, Fillard, Pennec and Ayache \cite{ArsignyLogE1,ArsignySIAM} and
Arsigny,  Pennec and Ayache \cite{ArsignyLogE2}.

\section{Appendix; Some Proofs Regarding the Jordan Form}
\label{sec6}
{\bf Proposition  \ref{elemdiv2}.}\  
{\it
Let $V$ be a finite-dimensional vector space 
and let $\mapdef{f}{V}{V}$ be a linear map.
If $V$ is a cyclic $\complex[X]$-module and if $(X - \lambda)^n$ is
the minimal polynomial of $f$, then there is a basis of $V$ of the form
\[
((f - \lambda \id)^{n-1}(u), (f - \lambda \id)^{n-2}(u), \ldots,
(f - \lambda \id)(u), u),
\]
for some $u\in V$. With respect to this basis, the matrix of $f$
is the Jordan block
\[
J_n(\lambda) = 
\begin{pmatrix}
\lambda & 1 & 0 &  \cdots &  0 \\
0 & \lambda  & 1  & \cdots & 0 \\
\vdots & \vdots  & \ddots  & \ddots  &  \vdots      \\
0 & 0 & 0 & \ddots &  1 \\
0 & 0 & 0 & \cdots &  \lambda
\end{pmatrix}.
\]
Consequently, $\lambda$ is an eigenvalue of $f$.
}

\begin{proof}
Since $V$ is a  cyclic $\complex[X]$-module, there is some $u\in V$ so that
$V$ is generated by  $u, f(u), f^2(u), \ldots$, 
which means that every vector in $V$ is of the form
$p(f)(u)$, for some polynomial, $p(X)$.
We claim that  $u, f(u),  \ldots, f^{n - 2}(u), f^{n-1}(u)$ generate $V$, 
which implies that the dimension of $V$ is at most $n$.

\medskip
This is because if $p(X)$ is any polynomial of degree
at least $n$, then we can divide $p(X)$ by $(X - \lambda)^n$ obtaining
\[
p = (X - \lambda)^n q + r,
\]
where $0 \leq \mathrm{deg}(r) < n$ and as $(X - \lambda)^n$ annihilates $V$,
we get
\[
p(f)(u) = r(f)(u),
\]
which means that every vector of the form $p(f)(u)$ with
$p(X)$ of degree $\geq n$ is actually a linear combination
of  $u, f(u),  \ldots, f^{n - 2}(u), f^{n-1}(u)$.

\medskip
We claim that the vectors
\[
u, (f - \lambda \id)(u),  \ldots, (f - \lambda \id)^{n-2}(u) 
(f - \lambda \id)^{n-1}(u)
\]
are linearly independent.
Indeed, if we had a nontrivial linear combination
\[
a_0(f - \lambda \id)^{n-1}(u) + a_1(f - \lambda \id)^{n-2}(u) + \cdots +
a_{n - 2}(f - \lambda \id)(u) + a_{n-1} u = 0, 
\]
then the polynomial
\[
a_0(X - \lambda)^{n-1} + a_1(X - \lambda)^{n-2} + \cdots +
a_{n - 2}(X - \lambda) + a_{n-1}   
\]
of degree at most $n -1$ would annihilate $V$, contradicting
the fact that $(X - \lambda)^n$ is the minimal polynomial of $f$
(and thus, of smallest degree).
Consequently, as the dimension of $V$ is at most $n$, 
\[
((f - \lambda \id)^{n-1}(u), (f - \lambda \id)^{n-2}(u), \ldots,
(f - \lambda \id)(u), u),
\]
is a basis of $V$ 
and since $u, f(u), \ldots, f^{n - 2}(u),  f^{n-1}(u)$ span $V$,
\[
(u, f(u),\ldots,  f^{n - 2}(u), f^{n-1}(u))
\]
is also a basis of $V$.  

\medskip
Let us see how $f$ acts on the basis
\[
((f - \lambda \id)^{n-1}(u), (f - \lambda \id)^{n-2}(u), \ldots,
(f - \lambda \id)(u), u).
\]

\medskip
If we write $f = f - \lambda\id + \lambda\id$, as
$(f - \lambda \id)^{n}$ annihilates $V$, we get
\[
f((f - \lambda \id)^{n-1}(u)) = (f - \lambda \id)^{n}(u) + 
\lambda (f - \lambda \id)^{n-1}(u) = \lambda (f - \lambda \id)^{n-1}(u)
\]
and
\[
f((f - \lambda \id)^{k}(u)) = (f - \lambda \id)^{k+1}(u) + 
\lambda (f - \lambda \id)^{k}(u),
\qquad 0 \leq k \leq n - 2.
\]
But this means precisely that the matrix of $f$ in this basis is the Jordan
block $J_n(\lambda)$.
\end{proof}

\medskip
To the best of our knowledge, a complete proof of the real Jordan form
is not easily found.
Horn and Johnson state such a result as
Theorem 3.4.5 in Chapter 3, Section 4, in \cite{HornJohn}.
However, they leave the details of the proof that a real $P$ can be found as
an exercise.  A complete proof is  given in Hirsh and Smale
\cite{HirshSmale}.  
This proof is given in Chapter 6,  and relies on results from Chapter
2 and  Appendix III. 

\medskip
We found that a proof  can be obtained from Theorem \ref{elemdiv1}. 
Since we believe that some of the techniques involved in this proof
are of independent interest, we present this proof in full detail.
It should be noted that we were inspired by some arguments found
in Gantmacher \cite{Gantmacher1}, Chapter IX, Section 13.

\medskip\noindent
{\bf Theorem \ref{Jordanform2}.} \ 
{\it 
(Real Jordan Form)
For every real $n\times n$ matrix, $A$, there is some
invertible (real) matrix, $P$, and some real Jordan matrix, $J$, so that
\[
A = P J P^{-1}.
\]
For every Jordan block, $J_r(\lambda)$, of type (1),
$\lambda$ is some real eigenvalue of $A$ and for
every Jordan block, $J_{2r}(\lambda, \mu)$, of type (2),
$\lambda + i \mu$ is a complex eigenvalue of $A$
(with $\mu\not= 0$). Every eigenvalue of $A$ corresponds to
one of more Jordan blocks of $J$.
Furthermore, the number, $m$, of Jordan blocks, 
the distinct Jordan block, $J_{s_i}(\alpha_i)$,  and
the number of times, $m_i$, that each Jordan block, $J_{s_i}(\alpha_i)$, occurs
are uniquely determined by $A$.
}

\begin{proof}
Let $\mapdef{f}{V}{V}$ be the linear map defined by $A$ and let $f_{\complex}$
be the complexification of $f$. Then, Theorem \ref{elemdiv1} yields
a direct sum decomposition of $V_{\complex}$ of the form
\begin{equation}
V_{\complex} = V_1 \oplus \cdots \oplus V_m,
\tag{$*$}
\end{equation}
where each $V_i$ is a cyclic $\complex[X]$-module 
(associated with $f_{\complex}$) whose minimal
polynomial is of the form $(X - \alpha_i)^{r_i}$,
where $\alpha$ is some (possibly complex) eigenvalue of $f$.
If $W$ is any subspace of $V_{\complex}$, we define the {\it conjugate\/},
$\overline{W}$,  of $W$ by
\[
\overline{W} = \{u - i v\in V_{\complex} \mid u + i v \in W\}.
\]
It is clear that $\overline{W}$ is a subspace of $V_{\complex}$ of the same 
dimension as $W$ and obviously, 
$\overline{V_{\complex}} = V_{\complex}$. 
Our first goal is to prove the following claim:

\medskip
{\it Claim 1\/}.
For each factor, $V_j$, the following properties hold:
\begin{enumerate}
\item[(1)]
If $u + i v, f_{\complex}(u + iv), \ldots, f_{\complex}^{r_j - 1} (u + iv)$
span $V_j$, then 
$u - i v, f_{\complex}(u - iv), \ldots, f_{\complex}^{r_j - 1} (u - iv)$
span $\overline{V}_i$ and so,
$\overline{V}_i$ is cyclic with respect to $f_{\complex}$.
\item[(2)]
If $(X - \alpha_i)^{r_i}$ is the minimal
polynomial of $V_i$, then $(X - \overline{\alpha}_i)^{r_i}$ is the minimal
polynomial of  $\overline{V}_i$. 
\end{enumerate}

\begin{proof}[Proof of Claim 1]
As $f_{\complex}(u + i v) = f(u) + i f(v)$, we have
$f_{\complex}(u - i v) = f(u) - i f(v)$. It follows that
$f_{\complex}^k(u + i v) = f^k(u) + i f^k(v)$ and
$f_{\complex}^k(u - i v) = f^k(u) - i f^k(v)$, which implies that
if $V_j$  is generated by
$u + i v, f_{\complex}(u + i v), \ldots, f_{\complex}^{r_j}(u + i v)$
then $\overline{V}_j$  is generated by \\
$u - i v, f_{\complex}(u - iv), \ldots, f_{\complex}^{r_j}(u - i v)$.
Therefore, $\overline{V}_j$ is cyclic for  $f_{\complex}$.

\medskip
We also prove the following simple fact:
If 
\[
(f_{\complex} - (\lambda_j + i\mu_j)\id)(u + iv) = x + i y, 
\]
then
\[
(f_{\complex} - (\lambda_j - i\mu_j)\id)(u - i v) = 
x - i y.
\]
Indeed, we have
\begin{eqnarray*}
 x + i y & = & (f_{\complex} - (\lambda_j + i\mu_j)\id)(u + iv) \\
 & = & f_{\complex}(u + iv)  - (\lambda_j + i\mu_j)(u + iv) \\
 & = & f(u) + i f(v)  - (\lambda_j + i\mu_j)(u + iv) 
\end{eqnarray*}
and by taking conjugates, we get
\begin{eqnarray*}
x - i y & = & f(u)  - if(v)  - (\lambda_j - i\mu_j)(u - iv)\\  
 & = & f_{\complex}(u - iv)  - (\lambda_j- i\mu_j)(u - iv) \\
 & = & (f_{\complex} - (\lambda_j - i\mu_j)\id)(u - iv),
\end{eqnarray*}
as claimed. 

\medskip
From the above,
$(f_{\complex} - \alpha_j\id)^{r_j}(x + iy) = 0$ iff
$(f_{\complex} - \overline{\alpha}_j\id)^{r_j}(x - iy) = 0$. 
Thus,  $(X -  \overline{\alpha}_j\id)^{r_j}$ annihilates  $\overline{V}_j$
and as $\mathrm{dim}\, \overline{V}_j = \mathrm{dim}\, V_j$
and  $\overline{V}_j$ is cyclic, we conclude that
$(X - \overline{\alpha}_j)^{r_j}$
is the minimal polynomial of $\overline{V}_j$. 
\end{proof}

\medskip
Next we prove

\medskip
{\it Claim 2\/}.
For every factor, $V_j$, in the direct decomposition $(*)$, we have:

\medskip
(A) If $(X - \lambda_j)^{r_j}$ is the minimal polynomial of 
$V_j$, with $\lambda_j \in \reals$, then either
\begin{enumerate}
\item[(1)]
$V_j = \overline{V}_j$ and if $u + i v$ generates $V_j$, then
$u - i v$ also generates $V_j$, or
\item[(2)]
$V_j \cap \overline{V}_j = (0)$ and
\begin{enumerate}
\item[(a)]
the cyclic space $\overline{V}_j$
also occurs in the direct sum decomposition $(*)$ 
\item[(b)]
the  minimal polynomial of  $\overline{V}_j$ is $(X - \lambda_j)^{r_j}$
\item[(c)]
the spaces $V_j$ and $\overline{V}_j$ contain
only complex vectors (this means that if $x + i y\in V_j$, then
$x\not= 0$ and $y\not= 0$ and similarly for  $\overline{V}_j$).
\end{enumerate}
\end{enumerate}

(B) If  $(X - (\lambda_j + i\mu_j))^{r_j}$ is the minimal polynomial of 
$V_j$ with $\mu_j \not= 0$, then
\begin{enumerate}
\item[(d)]
$V_j \cap \overline{V}_j = (0)$ 
\item[(e)]
the cyclic space $\overline{V}_j$
also occurs in the direct sum decomposition $(*)$
\item[(f)]
the  minimal polynomial of  $\overline{V}_j$ is 
$(X - (\lambda_j - i\mu_j))^{r_j}$
\item[(g)]
the spaces $V_j$ and $\overline{V}_j$ contain
only complex vectors.
\end{enumerate}

\begin{proof}[Proof of Claim 2]
By taking the conjugate
of the direct sum decomposition $(*)$ we get
\[
V_{\complex} = \overline{V}_1 \oplus \cdots \oplus \overline{V}_m.
\]
By Claim 1, each $\overline{V}_j$ is a cyclic subspace with 
respect to $f_{\complex}$ of the same dimension
as $V_j$ and the minimal polynomial of $\overline{V}_j$ is 
$(X - \overline{\alpha}_j)^{r_j}$ if the minimal polynomial of $V_j$ is
$(X - \alpha_j)^{r_j}$.
It follows from the uniqueness assertion of Theorem \ref{elemdiv1}
that the list  of conjugate minimal polynomials
\[
(X - \overline{\alpha}_1)^{r_1}, \ldots, (X - \overline{\alpha}_m)^{r_m}
\]
is a permutation the list of minimal polynomials
\[
(X - \alpha_1)^{r_1}, \ldots, (X - \alpha_m)^{r_m}
\]
and so,  every $\overline{V}_j$ is equal to some factor $V_k$
(possibly equal to $V_j$ if $\alpha_j$ is real)
in the direct decomposition $(*)$, where $V_k$ and $\overline{V}_j$
have the same minimal polynomial, $(X - \overline{\alpha}_j)^{r_j}$.

\medskip
Next, assume that $(X - \lambda_j)^{r_j}$ is the minimal polynomial of 
$V_j$, with $\lambda_j \in \reals$. Consider any generator,
$u + i v$, for $V_j$. If $u - i v\in V_j$, then by Claim 1,
$\overline{V}_j \subseteq  V_j$  and so $\overline{V}_j =  V_j$, 
as  $\mathrm{dim}\, \overline{V}_j = \mathrm{dim}\, V_j$.
We know that 
$u + i v, f_{\complex}(u + i v), \ldots, f_{\complex}^{r_j}(u + i v)$ 
generate $V_j$ and that
$u - i v, f_{\complex}(u - iv), \ldots, f_{\complex}^{r_j}(u - i v)$
generate $\overline{V}_j = V_j$, which implies (1).

\medskip
If $u - iv\notin V_j$, then we proved earlier that $\overline{V}_j$
occurs in the direct sum $(*)$ as some $V_k$ 
and that its minimal polynomial is also  $(X - \lambda_j)^{r_j}$.
Since $u - iv\notin V_j$ and $V_j$ and $\overline{V}_j$
belong to a direct sum decomposition,  $V_j \cap \overline{V}_j = (0)$
and 2(a) and 2(b) hold.
If $u\in V_j$  or $iv\in V_j$ for some real $u\in V$ or some real $v\in V$
and $u, v\not= 0$, 
as $V_j$ is a complex space, then $v\in V_j$ and
either $u \in \overline{V}_j$ or $v \in \overline{V}_j$,
contradicting  $V_j \cap \overline{V}_j = (0)$. Thus, 2(c) holds.

\medskip
Now, consider the case where $\alpha_j = \lambda_j + i\mu_j$,
with $\mu_j\not= 0$. Then, we know that
$\overline{V}_j = V_k$ for some
$V_k$ whose  minimal polynomial is 
$(X - (\alpha_j - i\mu_j))^{r_j}$ in the direct sum $(*)$.
As $\mu_j\not= 0$, the cyclic spaces $V_j$ and $\overline{V}_j$
correspond to distinct minimal polynomials  $(X - (\alpha_j + i\mu_j))^{r_j}$
and  $(X - (\alpha_j - i\mu_j))^{r_j}$, so
$V_j\cap \overline{V}_j = (0)$. It follows that $V_j$ and $\overline{V}_j$
consist of complex vectors as we already observed.
Therefore, (d), (e), (f), (g) are proved, which finishes the proof
of Claim 2.
\end{proof}
This completes the proof our theorem.
\end{proof}

\medskip\noindent
{\bf Theorem \ref{expJord}.} \   
{\it 
For any (real or complex) $n\times n$ matrix, $A$, if
$A = P J P^{-1}$ where $J$ is a Jordan matrix of the form
\[
J = 
\begin{pmatrix}
J_{r_1}(\lambda_1) &  \cdots & 0 \\
 \vdots  & \ddots   & \vdots    \\
 0 & \cdots & J_{r_m}(\lambda_m)
\end{pmatrix},
\]
then there is some invertible matrix, $Q$, so that
the  Jordan form of $e^A$ is given by
\[
e^A = Q\, e(J)\,  Q^{-1},
\]
where $e(J)$ is the Jordan matrix
\[
e(J) = 
\begin{pmatrix}
J_{r_1}(e^{\lambda_1}) &  \cdots & 0 \\
 \vdots  & \ddots   & \vdots    \\
 0 & \cdots & J_{r_m}(e^{\lambda_m})
\end{pmatrix},
\]
that is, each $J_{r_k}(e^{\lambda_k})$ is obtained from $J_{r_k}(\lambda_k)$
by replacing all the diagonal entries $\lambda_k$ by $e^{\lambda_k}$.
Equivalently, if the list of elementary divisors of $A$ is
\[
(X - \lambda_1)^{r_1}, \ldots, (X - \lambda_m)^{r_m},
\]
then the list of elementary divisors of $e^A$ is
\[
(X - e^{\lambda_1})^{r_1}, \ldots, (X - e^{\lambda_m})^{r_m}.
\]
}

\begin{proof}
Theorem \ref{expJord} is a consequence of a general theorem
about functions of matrices proved in Gantmacher \cite{Gantmacher1},
see Chapter VI, Section 8, Theorem 9. Because a much more 
general result is proved, the proof in Gantmacher \cite{Gantmacher1}
is rather involved.
However, it is possible to give a simpler proof exploiting
special properties of the exponential map.

\medskip
Let $f$ be the linear map defined by the matrix $A$. 
The strategy of our proof is to go back to the direct sum decomposition
given by Theorem \ref{elemdiv1},
\[
V = V_1 \oplus V_2 \oplus \cdots \oplus V_m,
\]
where each $V_i$ is a cyclic $\complex[X]$-module 
such that the minimal polynomial
of the restriction of $f$ to $V_i$ is of the form $(X - \lambda_i)^{r_i}$.
We will prove that
\begin{enumerate}
\item[(1)]
The vectors
\[
u, e^f(u), (e^f)^2(u), \ldots, (e^f)^{r_i - 1}(u)
\]
form a basis of $V_i$ (here, $(e^f)^k = e^f \circ \cdots \circ e^f$, the
composition of $e^f$ with itself $k$ times).
\item[(2)]
The polynomial
$(X - e^{\lambda_i})^{r_i}$ is
the minimal polynomial of the restriction of $e^f$ to $V_i$.
\end{enumerate}

First, we prove that $V_i$ is invariant under $e^f$.
Let $N = f - \lambda_i \id$. To say that $(X - \lambda_i)^{r_i}$
is the minimal polynomial of the restriction of $f$ to $V_i$ is equivalent to
saying that $N$ is nilpotent with index of nilpotency, $r = r_j$.
Now, $N$  and $\lambda_i \id$ commute so as $f = N + \lambda_i \id$,
we have
\[
e^f = e^{N + \lambda_i \id} = e^{N} e^{\lambda_i \id} = e^{\lambda_i} e^{N}. 
\]
Furthermore, as $N$ is nilpotent, we have
\[
e^N = \id + N + \frac{N^2}{2!} + \cdots + \frac{N^{r-1}}{(r -1)!},
\]
so
\[
e^f =  e^{\lambda_i}
\left(
\id + N + \frac{N^2}{2!} + \cdots + \frac{N^{r-1}}{(r -1)!}
\right).
\]
Now, $V_i$ is invariant under $f$ so $V_i$ is
invariant under $N =  f - \lambda_i \id$
and this implies that $V_i$ is invariant under $e^f$.
Thus, we can view $V_i$ as a $\complex[X]$-module with respect to $e^f$.

\medskip
From the formula for $e^f$ we get
\begin{eqnarray*}
e^f - e^{\lambda_i} \id & = & e^{\lambda_i}
\left(
\id + N + \frac{N^2}{2!} + \cdots + \frac{N^{r-1}}{(r -1)!}
\right) - e^{\lambda_i}\id \\
& = & e^{\lambda_i}
\left(
 N + \frac{N^2}{2!} + \cdots + \frac{N^{r-1}}{(r -1)!}
\right).
\end{eqnarray*}
If we let 
\[
\widetilde{N} =  N + \frac{N^2}{2!} + \cdots + \frac{N^{r-1}}{(r - 1)!},
\]
we claim that
\[
\widetilde{N}^{r - 1} = N^{r -1}
\quad\hbox{and}\quad
\widetilde{N}^{r} = 0. 
\]

\medskip
The case $r = 1$ is trivial so we may assume $r \geq 2$.
Since  $\widetilde{N} = N R$ for some $R$ such that
$N R = R N$ and $N^r = 0$,
the second property is clear. The first property
follows by observing 
that $\widetilde{N} = N + N^2T$, where $N$ and $T$ commute, so using
the binomial formula,
\[
\widetilde{N}^{r - 1} = 
\sum_{k = 0}^{r - 1} \binom{r-1}{k} N^k (N^2T)^{r - 1 - k}
= \sum_{k = 0}^{r - 1} \binom{r-1}{k} N^{2r - k - 2}T^{r - 1 - k} = N^{r-1},
\]
since $2r - k - 2 \geq r$ for $0 \leq k \leq r - 2$ and $N^r = 0$.

\medskip
Recall from Proposition \ref{elemdiv2} that 
\[
((f - \lambda_i\id)^{r_i - 1}(u), \ldots, (f - \lambda_i\id)(u), u)
\]
is a basis of $V_i$, which implies that
$N^{r - 1}(u) = (f - \lambda_i\id)^{r_i - 1}(u) \not= 0$. 
Since  $\widetilde{N}^{r - 1} = N^{r -1}$,
we have  $\widetilde{N}^{r - 1}(u) \not= 0$ and as
$\widetilde{N}^{r} = 0$, we have $\widetilde{N}^{r}(u)  = 0$.
It is well-known that these two facts imply that
\[
u, \widetilde{N}(u), \ldots, \widetilde{N}^{r-1}(u)
\]
are linearly independent. Indeed, if we had a linear dependence relation
\[
a_0 u + a_1\widetilde{N}(u) +  \cdots +  a_{r-1}\widetilde{N}^{r-1}(u) = 0,
\]
by applying $\widetilde{N}^{r - 1}$, as  $\widetilde{N}^{r}(u)  = 0$
we get $a_0\widetilde{N}^{r - 1}(u) = 0$,
so, $a_0 = 0$ as $\widetilde{N}^{r - 1}(u) \not= 0$; by 
applying $\widetilde{N}^{r - 2}$ we get  $a_1\widetilde{N}^{r - 1}(u) = 0$,
so $a_1 = 0$;  using induction, by applying  $\widetilde{N}^{r - k - 2}$ to
\[
a_{k+1}\widetilde{N}^{k+1}(u) +  \cdots +  a_{r - 1}\widetilde{N}^{r-1}(u) = 0,
\]
we get $a_{k+1}= 0$ for $k = 0, \ldots, r - 2$.
Since $V_i$ has dimension $r\> (= r_i)$, we deduce that
\[
(u, \widetilde{N}(u), \ldots, \widetilde{N}^{r-1}(u))
\]
is a basis of $V_i$. But $e^f = e^{\lambda_i}(\id + \widetilde{N})$,
so for $k = 0, \ldots, r-1$, each
$\widetilde{N}^{k}(u)$ is a linear combination of
the vectors $u, e^f(u), \ldots, (e^f)^{r- 1}(u)$ which implies that
\[
(u, e^f(u), (e^f)^2(u), \ldots, (e^f)^{r- 1}(u))
\]
is a basis of $V_i$. This implies that any annihilating
polynomial of $V_i$ has degree no less than $r$ and since
$(X - e^{\lambda_i})^r$ annihilates $V_i$ (because
$e^f - e^{\lambda_i}\id =  e^{\lambda_i}\widetilde{N}$ and 
$\widetilde{N}^r = 0$), it is the minimal
polynomial of $V_i$. 

\medskip
In summary, we proved that each $V_i$ is a cyclic $\complex[X]$-module
(with respect to $e^f$) and that in the direct sum decomposition
\[
V = V_1\oplus \cdots \oplus V_m,
\]
the polynomial $(X - e^{\lambda_i})^{r_i}$ is the minimal polynomial
of $V_i$, which is Theorem  \ref{elemdiv1} for $e^f$.
Then, Theorem \ref{expJord} follows immediately from
Proposition \ref{elemdiv2}.
\end{proof}

\medskip\noindent
{\bf Theorem \ \ref{squareJord}.} \ 
{\it 
For any (real or complex) invertible $n\times n$ matrix, $A$, if
$A = P J P^{-1}$ where $J$ is a Jordan matrix of the form
\[
J = 
\begin{pmatrix}
J_{r_1}(\lambda_1) &  \cdots & 0 \\
 \vdots  & \ddots   & \vdots    \\
 0 & \cdots & J_{r_m}(\lambda_m)
\end{pmatrix},
\]
then there is some invertible matrix, $Q$, so that
the  Jordan form of $A^2$ is given by
\[
e^A = Q\, s(J)\,  Q^{-1},
\]
where $s(J)$ is the Jordan matrix
\[
s(J) = 
\begin{pmatrix}
J_{r_1}(\lambda_1^2) &  \cdots & 0 \\
 \vdots  & \ddots   & \vdots    \\
 0 & \cdots & J_{r_m}(\lambda_m^2)
\end{pmatrix},
\]
that is, each $J_{r_k}(\lambda_k^2)$ is obtained from $J_{r_k}(\lambda_k)$
by replacing all the diagonal enties $\lambda_k$ by $\lambda_k^2$.
Equivalently, if the list of elementary divisors of $A$ is
\[
(X - \lambda_1)^{r_1}, \ldots, (X - \lambda_m)^{r_m},
\]
then the list of elementary divisors of $A^2$ is
\[
(X - \lambda_1^2)^{r_1}, \ldots, (X - \lambda_m^2)^{r_m}.
\]
}

\begin{proof}
Theorem \ref{squareJord} is a consequence of a general theorem
about functions of matrices proved in Gantmacher \cite{Gantmacher1},
see Chapter VI, Section 8, Theorem 9. 
However, it is possible to give a simpler proof exploiting
special properties of the squaring map.

\medskip
Let $f$ be the linear map defined by the matrix $A$.
The proof is modeled after the proof of Theorem  \ref{expJord}.
Consider the direct sum decomposition
given by Theorem \ref{elemdiv1},
\[
V = V_1 \oplus V_2 \oplus \cdots \oplus V_m,
\]
where each $V_i$ is a cyclic $\complex[X]$-module 
such that the minimal polynomial
of the restriction of $f$ to $V_i$ is of the form $(X - \lambda_i)^{r_i}$.
We will prove that
\begin{enumerate}
\item[(1)]
The vectors
\[
u, f^2(u), f^4(u), \ldots, f^{2(r_i - 1)}(u)
\]
form a basis of $V_i$.
\item[(2)]
The polynomial
$(X - \lambda_i^2)^{r_i}$ is
the minimal polynomial of the restriction of $f^2$ to $V_i$.
\end{enumerate}

Since $V_i$ is invariant under $f$, it is clear that $V_i$ is invariant
under $f^2 = f\circ f$.
Thus, we can view $V_i$ as a $\complex[X]$-module with respect to $f^2$.
Let $N = f - \lambda_i \id$. To say that $(X - \lambda_i)^{r_i}$
is the minimal polynomial of the restriction of $f$ to $V_i$ is equivalent to
saying that $N$ is nilpotent with index of nilpotency, $r = r_j$.
Now, $N$  and $\lambda_i \id$ commute so as $f = \lambda_i \id + N$,
we have
\[
f^2 = \lambda_i^2 \id + 2\lambda_i N + N^2
\]
and so
\[
f^2 - \lambda_i^2 \id =  2\lambda_i N + N^2.
\]
Since we are assuming that $f$ is invertible, $\lambda_i \not= 0$, so 
\[
f^2 - \lambda_i^2 \id =  2\lambda_i \left(N + \frac{N^2}{2\lambda_i}\right).
\]
If we let
\[
\widetilde{N} =  N + \frac{N^2}{2\lambda_i},
\]
we claim that
\[
\widetilde{N}^{r - 1} = N^{r -1}
\quad\hbox{and}\quad
\widetilde{N}^{r} = 0. 
\]

\medskip
The proof is identical to the proof given in Theorem \ref{expJord}.
Again, as in the proof of  Theorem \ref{expJord}, we deduce that
we have  $\widetilde{N}^{r - 1}(u) \not= 0$ and $\widetilde{N}^{r}(u)  = 0$, from
which we infer that
\[
(u, \widetilde{N}(u), \ldots, \widetilde{N}^{r-1}(u))
\]
is a basis of $V_i$. But $f^2 - \lambda_i^2 \id = 2\lambda_i\widetilde{N}$,
so for $k = 0, \ldots, r-1$, each
$\widetilde{N}^{k}(u)$ is a linear combination of
the vectors $u, f^2(u), \ldots, f^{2(r- 1)}(0)$ which implies that
\[
(u, f^2(u), f^4(u), \ldots, f^{2(r- 1)}(u))
\]
is a basis of $V_i$. This implies that any annihilating
polynomial of $V_i$ has degree no less than $r$ and since
$(X - \lambda_i^2)^r$ annihilates $V_i$
(because  $f^2 - \lambda_i^2 \id = 2\lambda_i\widetilde{N}$ and
$\widetilde{N}^r = 0$), it is the minimal
polynomial of $V_i$. 
Theorem \ref{squareJord} follows immediately from
Proposition \ref{elemdiv2}.
\end{proof}


\bibliographystyle{plain} 

\end{document}